\documentclass[authordraft=false, format=acmsmall, review=false, screen=true, nonacm=true]{acmart}

\AtBeginDocument{%
  \providecommand\BibTeX{{%
      \normalfont B\kern-0.5em{\scshape i\kern-0.25em b}\kern-0.8em\TeX}}}

\usepackage{amsmath}
\usepackage{amsfonts}

\usepackage[binary-units=true]{siunitx}
\usepackage{subcaption} % subfigure, \ContinuedFloat
\usepackage{enumitem} % [resume]
\usepackage{isomath} % \vectorsym, \matrixsym, \tensorsym

\usepackage{glossaries-extra} % \newabbreviation, gls
\glsdisablehyper

\usepackage{adjustbox} % for colorbar

\usepackage{femtikz} % pgfplots, tikz, graphicx, \LagrangeCell
\pgfplotsset{compat=1.9} % colorbar style ylabel
\usepgfplotslibrary{colorbrewer, colormaps, groupplots}
\pgfplotsset{
  cycle list/Dark2,
  cycle multiindex* list={
    mark list*\nextlist
    Dark2\nextlist
  },
  every axis plot/.append style={thick}
}

\newcommand{\ghostfont}{\itshape}

\newcommand{\h}{\textit{h}}
\newcommand{\p}{\textit{p}}
\newcommand{\hp}{\textit{hp}}

\newcommand{\gmres}{GMRES}
\newcommand{\expanse}{Expanse}

\newcommand{\dealii}{\texttt{deal.II}}
\newcommand{\hpbox}{\texttt{hpbox}}

\newcommand{\pforest}{\texttt{p4est}}

\newcommand{\trilinos}{\texttt{Trilinos}}

\newcommand{\ml}{\texttt{ML}}

\newcommand{\dune}{\texttt{DUNE}}
\newcommand{\fenics}{\texttt{FEniCS}}
\newcommand{\fenicsx}{\texttt{FEniCSx}}
\newcommand{\freefem}{\texttt{FreeFEM++}}
\newcommand{\getfem}{\texttt{GetFEM}}
\newcommand{\hermes}{\texttt{Hermes}}
\newcommand{\libmesh}{\texttt{libMesh}}
\newcommand{\mfem}{\texttt{MFEM}}
\newcommand{\mofem}{\texttt{MoFEM}}
\newcommand{\phg}{\texttt{PHG}}
\newcommand{\phaml}{\texttt{PHAML}}

\newtheorem{remark}{Remark}

\renewcommand{\vec}[1]{\ensuremath{\vectorsym{#1}}}

\newcommand{\cellsp}[1]{\ensuremath{\mathbb{T}^p_{\text{#1}}}}

\newabbreviation{cpu}{CPU}{central processing unit}
\newabbreviation{gpu}{GPU}{graphics processing unit}

\newabbreviation{api}{API}{application programming interface}
\newabbreviation{cuda}{CUDA}{Compute Unified Device Architecture}
\newabbreviation{mpi}{MPI}{Message Passing Interface}
\newabbreviation{openacc}{OpenACC}{Open Accelerators}
\newabbreviation{openmp}{OpenMP}{Open Multi-Processing}
\newabbreviation{tbb}{TBB}{Threading Building Blocks}

\newabbreviation{cfd}{CFD}{computational fluid dynamics}

\newabbreviation{hpc}{HPC}{high-performance computing}

\newabbreviation[plural={FDM}, \glsshortpluralkey={FDM}]{fdm}{FDM}{finite difference method}
\newabbreviation[plural={FVM}, \glsshortpluralkey={FVM}]{fvm}{FVM}{finite volume method}
\newabbreviation[plural={FEM}, \glsshortpluralkey={FEM}]{fem}{FEM}{finite element method}
\newabbreviation{fft}{FFT}{fast Fourier transformation}
\newabbreviation{lbm}{LBM}{lattice Boltzmann method}

\newabbreviation[longplural={degrees of freedom}]{dof}{DoF}{degree of freedom}

\newabbreviation{amr}{AMR}{adaptive mesh refinement}
\newabbreviation{simd}{SIMD}{single instruction, multiple data}

\newabbreviation{cg}{CG}{continuous Galerkin}
\newabbreviation{dg}{DG}{discontinuous Galerkin}

\newabbreviation{amg}{AMG}{algebraic multigrid}
\newabbreviation{gmg}{GMG}{geometric multigrid}

\newabbreviation{csr}{CSR}{compressed row storage}

\newabbreviation{xsdk}{xSDK}{extreme-scale scientific software development kit}
\begin{document}
% Title portion
\title{Algorithms for Parallel Generic \hp-adaptive Finite Element Software}
\titlenote{Dedicated to the memory of William F. Mitchell.}
% \subtitle{This is a subtitle}
% \subtitlenote{Subtitle note}

\author{Marc Fehling}
\orcid{0000-0003-0984-793X}
\affiliation{%
  \department{Department of Mathematics}
  \institution{Colorado State University}
  \streetaddress{1874 Campus Delivery}
  \city{Fort Collins}
  \state{CO}
  \postcode{80523-1874}
  \country{USA}}
\email{marc.fehling@colostate.edu}

\author{Wolfgang Bangerth}
\orcid{0000-0003-2311-9402}
\affiliation{%
  \department{Department of Mathematics and Department of Geosciences}
  \institution{Colorado State University}
  \streetaddress{1874 Campus Delivery}
  \city{Fort Collins}
  \state{CO}
  \postcode{80523-1874}
  \country{USA}}
\email{bangerth@colostate.edu}

\renewcommand\shortauthors{M. Fehling and W. Bangerth}

\begin{abstract}
The \hp-adaptive \gls{fem} -- where one independently chooses the mesh size (\h) and polynomial degree (\p) to be used on each cell -- has long been known to have better theoretical convergence properties than either \h- or \p-adaptive methods alone. However, it is not widely used, owing at least in parts to the difficulty of the underlying algorithms and the lack of widely usable implementations. This is particularly true when used with continuous finite elements.

Herein, we discuss algorithms that are necessary for a comprehensive and generic implementation of \hp-adaptive finite element methods on distributed-memory, parallel machines. In particular, we will present a multi-stage algorithm for the unique enumeration of \glspl{dof} suitable for continuous finite element spaces, describe considerations for weighted load balancing, and discuss the transfer of variable size data between processes. We illustrate the performance of our algorithms with numerical examples, and demonstrate that they scale reasonably up to at least \num{16384} \gls{mpi} processes.

We provide a reference implementation of our algorithms as part of the open-source library \dealii{}.
\end{abstract}






\maketitle

% Content
\glsresetall
\section{Introduction}

In the \hp-adaptive variation of the \gls{fem} for the solution of partial differential equations, one adaptively refines the mesh (\h-adaptivity) and independently also chooses the polynomial degree of the approximation on every cell (\p-adaptivity). This method is by now 40 years old \cite{babuska1981} and, at least from a theoretical perspective, well understood \cite{guo1986,guo1986a,babuska1996}. In particular, it is known that \hp-adaptivity provides better accuracy per \gls{dof} than either the \h- or \p-adaptive methods alone; more specifically, it exhibits a convergence rate where the approximation error in many cases decreases \textit{exponentially} with the number of unknowns $N$ -- i.e., the error satisfies $e = {\mathcal O}(s^{-N})$ for some $s>1$ that may depend on the solution --, rather than an algebraic rate $e = {\mathcal O}(N^{-\gamma})$ for some $\gamma>0$. In other words, \hp-adaptivity is \textit{asymptotically superior} to \h- or \p-adaptivity alone.

Yet, \textit{\hp-adaptive methods are not widely used}. The reasons for this lack of use are probably debatable but surely include (i) that the literature provides many criteria by which to choose whether \h- or \p-refinement should be selected if the error on a cell is large, but that there is no consensus on which one is best; and (ii) a lack of widely usable implementations. For the first of these points, we refer to the comprehensive comparison in \cite{mitchell2014}. Instead, in this contribution, we address the second point: the lack of widely available implementations.

A survey of the finite element landscape shows that there are few options for those who are interested in experimenting with \hp-methods. Most of the open-source distributed-memory parallel implementations of \hp-adaptive methods available that we are aware of -- specifically the ones in the libraries \phaml{} \cite{mitchell2002}, \phg{} \cite{zhang2019}, and \mofem{} \cite{kaczmarczyk2020} -- have not found wide use in the community and are not backed by large user and developer communities. To the best of our knowledge, other popular libraries like \fenics{}/\fenicsx{} \cite{alnaes2015}, \getfem{} \cite{renard2020}, and \freefem{} \cite{hecht2012} do not offer \hp-adaptive methods at all or have only experimental support as is the case with \libmesh{} \cite{kirk2006}.

In other cases, such as the ones discussed in \cite{bey1996, paszynski2006, paszynski2011, chalmers2019}, the implementation of \hp-methods is restricted to \gls{dg} methods; the same limitation also applies to the libraries \mfem{} \cite{anderson2021, pazner2022} and \dune{} \cite{bastian2021, gersbacher2016}. This case is relatively easy to implement because the construction of finite element spaces is purely local, on every cell independent of its neighbors. At the same time, \gls{dg} methods are expensive -- especially in three dimensions -- because \glspl{dof} are duplicated between neighboring cells, and the resulting large linear systems and corresponding memory consumption have hampered adoption of \gls{dg} schemes in most applications outside the simulation of hyperbolic systems. As a consequence, while the use of \gls{dg} methods for \hp-adaptivity is a legitimate approach, there are many important use cases where continuous finite element spaces remain the method of choice.

Finally, let us mention publications \cite{paszynski2011, jomo2017}
that also demonstrate the use of \hp-adaptive methods; these use
distributed-memory parallelization, but use data structures for the
mesh replicated on all processes, thus limiting scalability. An
extension to distributed data structures, using a
hierarchical construction of finite element spaces, is discussed in \cite{jomo2021phd}. The \hermes{} library \cite{solin2008} falls into a separate category, using shared-memory parallelism. \cite{laszloffy2000} does present distributed-memory algorithms, but only shows scaling to 16 processors, whereas we are here interested in much larger levels of parallelism. We are not aware of any commercial tools capable of using \hp-methods, either for sequential or parallel computations.

As a consequence of our search for available implementations, and to the best of our knowledge, only the \dealii{} library \cite{dealii94,arndt2021} appears to have generic support for \hp-adaptive methods for a wide variety of finite elements, discontinuous or continuous, as previously discussed in detail in \cite{bangerth2009}. Still, \dealii{} has only recently begun to support \hp-adaptive methods for parallel computations \cite{fehling2020}. It is this specific gap that we wish to address in this contribution, by considering what algorithms are necessary to implement \hp-methods \textit{on large parallel machines} using a distributed-memory model based on the \gls{mpi}. The target for our work is the solution of two- and three-dimensional scalar- or vector-valued partial differential equations, using an arbitrary combination of finite elements, and scaling up to tens of thousands of processes and billions of unknowns.

%More specifically, our goals for this work are:
More specifically, we identify and address the following three major challenges in this work:
\begin{itemize}
    \item The development of a scalable algorithm to \textit{uniquely enumerate \glspl{dof}} on meshes on which finite element spaces of different polynomial degrees may be associated with each cell. Simply enumerating all \glspl{dof} on a mesh turns out to be non-trivial already in distributed-memory implementations of \h-adapted, unstructured meshes (as discussed in \cite{bangerth2012}) as well as for sequential implementations of the \hp-method (see \cite{bangerth2007}), and it is no surprise that the combination of the two leads to additional complications.
    
    \item An efficient distribution of workload among all processes with \textit{weighted load balancing}, since the workload per cell depends on its local number of \glspl{dof} and thus varies from cell to cell with \hp-adaptive methods. We will present strategies on how to determine weights on each cell for this purpose.
    
    \item The ability to \textit{transfer data of variable size} between \hp-adapted meshes during repartitioning. In the \hp-context, the amount of data stored per cell is proportional to the number of local \glspl{dof} and, again, varies between cells.
    
    \item An assessment of the parallel efficiency of the algorithms mentioned above.
\end{itemize}

% On the other hand, our goals do not include a discussion of \textit{how} to decide which cells to adapt. Instead, we refer to \cite{mitchell2014} for this question.

In this paper, we will first address the task of enumerating all \glspl{dof} in a distributed-memory setting in Section~\ref{sec:enumeration}. We will then present strategies for weighted load balancing in Section~\ref{sec:load_balancing} and continue with ways to transfer data of variable size in Section~\ref{sec:data_transfer}. In Section~\ref{sec:results}, we then illustrate the performance and scalability of our methods using numerical results obtained on the \expanse{} supercomputer \cite{strande2021}, using up to \num{16384} cores. We present conclusions in Section~\ref{sec:conclusions}.

\paragraph*{Code availability.}
The algorithms we will discuss in the remainder of this paper are implemented and available in the open-source library \dealii{}, version 9.4 \cite{dealii94,arndt2021}. All functionality is available under the LGPL 2.1 license. That said, our discussions are not specific to \dealii{} and are generally applicable to any other finite element software. In particular, even though we will only show examples of quadrilateral or hexahedral meshes, our algorithms are readily applicable also to simplex or mixed meshes.

The two programs that implement the test cases of Section~\ref{sec:testcases} and for which we show results in Sections~\ref{sec:results-load-balancing} and \ref{sec:results-scaling} are available as part of the tool \hpbox{} \cite{fehling2022}.

% -----------------------------------------
% --- Enumeration of degrees of freedom ---
% -----------------------------------------

\section{Enumeration of degrees of freedom}
\label{sec:enumeration}

In the abstract, the finite element method defines a finite-dimensional space $V_h$ within which one seeks the discrete solution of a (partial) differential equation. In practice, one needs to construct a basis $\{\varphi_i\}_{i=0}^{N-1}$ for this space so that numerical solutions $u_h\in V_h$ can be expressed as expansions of the form $u_h(\mathbf x) = \sum_i U_i \varphi_i(\mathbf x)$ where the $U_i$ are the nodal coefficients of the expansion.

The basis functions of $V_h$ are mathematically defined via \textit{nodal functionals} \cite{brenner2008}, but for the purposes of this section, it is only important to know that each basis function is associated with either a vertex, an edge, a face, or the interior of a cell of a mesh. In order to enumerate the \glspl{dof} on an unstructured mesh, one therefore simply walks over all cells, faces, edges, and vertices and, in a first step, allocates as much memory as is necessary to store the indices of \glspl{dof} associated with each of these entities, setting the index to an invalid value. In a second step, one then repeats the loop and assigns consecutive indices to each degree encountered that has not yet received a valid number.

Our goal herein is to come up with an algorithm that replicates this action for parallel computations if the data structures that represent the mesh and the indices of \glspl{dof} defined thereon are distributed across individual nodes of a parallel computer. Implicit in this goal is that we continue to want a single, global mesh along with a \textit{global enumeration} of all \glspl{dof} -- even though each process in this distributed memory universe only sees its own small part of this distributed data structure. There are of course other ways to solve partial differential equations in parallel, most notably using domain decomposition methods in which each process has only a local enumeration of the \glspl{dof} located on that part of the mesh it owns (and potentially on one or more layers of ``ghost cells''), plus a way to map its local indices to the local indices on neighboring processes at partition interfaces and on ghost cells. Yet, domain decomposition methods have somewhat fallen out of favor because they do not scale well to very large process counts, and have largely been replaced by methods that instead consider a global finite element space (with basis functions indexed globally) from which one then builds a single, global linear system that is stored in a distributed data structure spread across processes. It is this latter model we seek to support in our work, requiring a global enumeration of all \glspl{dof} with globally unique identifiers.

In the remainder of this section, our goal is to describe an algorithm that achieves this enumeration of \glspl{dof} in parallel for the \hp-adaptive case. For context, let us first briefly outline how this is done for distributed, unstructured meshes when only one type of finite element is used (Section~\ref{sec:enum-parallel}), followed by a description of the algorithm used for \hp-adaptive methods on a single process  (Section~\ref{sec:enum-hp}). In Section~\ref{sec:enum-hp-parallel}, we then present our new algorithm for parallel \hp-adaptive methods, which can be seen as a combination and enhancement of the former two.

We do not cover details on handling hanging nodes and constraints in this manuscript. It turns out that for the new algorithm, their handling does not require any change from the methods described in \cite{bangerth2009, bangerth2012}.

\subsection{Enumerating degrees of freedom on distributed, unstructured meshes}
\label{sec:enum-parallel}

In a parallel program where the mesh data structure is stored in distributed memory, the situation is complicated by the fact that each process only knows a subset of cells -- namely, those cells that are ``locally owned'' along with a layer of ``ghost cells''. At the same time, we need to assign \textit{globally unique} indices to all entities of the distributed mesh: at the end of the algorithm, each process must know the global indices of those \glspl{dof} that are located on this process's locally owned and ghost cells. 

For the relatively simple case where the finite element is the same on each cell (no \p-adaptivity), the index assignment is typically achieved by identifying a tie-breaking process that defines which process ``owns'' a mesh entity on the interface between the sub-domains of cells owned by individual processes (i.e., which of the adjacent processes owns a vertex, an edge, or a face on this interface). This process is then also the owner of the \glspl{dof} located on these entities. A possible tie-breaker is that the process with the smallest \gls{mpi} rank is chosen as the owner of an entity on a subdomain interface.

Enumeration of \glspl{dof} then proceeds by each process enumerating
the \glspl{dof} it owns, starting at zero. All of these indices are
then shifted so that we obtain globally unique indices across
processes. Next, each process sends the indices associated with
locally owned cells to those processes that have these cells as ghost
cells. Because processes may not yet know all \gls{dof} indices on the
boundaries of locally owned cells at the time of this communication step, the exchange has to be repeated a second time to ensure that each process knows the full set of indices on both the locally owned cells as well as ghost cells, and all of the vertices, edges, and faces bounding these cells.

A formal description of this algorithm -- which consists of five stages -- has been given in \cite{bangerth2012} and forms the basis of the discussions for the parallel \hp-case below in Section~\ref{sec:enum-hp-parallel}.

\subsection{Enumerating degrees of freedom in the sequential \hp-context}
\label{sec:enum-hp}

In the \hp-context, each cell $K\in {\mathbb T}$ of a triangulation or mesh $\mathbb T$ may use a different finite element. To make the notation that we use below concrete, let us assume that we want the global function space $V_h$ be constructed so that the solution functions $u_h\in V_h$ satisfy $u_h|_K \in V_h(K)$ where $V_h(K)$ is the finite element space associated with cell $K$. Furthermore, let us assume that $V_h(K)$ can only be one from within a collection of spaces $\bigl\{\hat V_h^{(i)}\bigr\}_{i=0}^I$ defined on the reference cell $\hat K$ that are then mapped to cell $K$ in the usual way, i.e., $V_h(K)={\mathcal M}_K \hat V_h$ where ${\mathcal M}_K$ is the operator that maps the finite element space from the reference cell to $K$; the details of this mapping are not of importance to us here. We denote the ``active FE index'' on cell $K$ by $a(K)$, i.e., $V_h(K)={\mathcal M}_K \hat V_h^{\left(a(K)\right)}$. Each of the spaces $\hat V_h^{(i)}$ has a number of \glspl{dof} associated with each vertex, edge, face, and cell interior.

A trivial implementation of enumerating all \glspl{dof} would simply loop over all cells $K\in\mathbb T$ and enumerate all \glspl{dof} on both the cell $K$ and its vertices, edges, and faces independently of the enumeration on neighboring cells. To do so requires storing multiple sets of indices of \glspl{dof} on vertices, edges, and faces, each set corresponding to one of the adjacent cells. This strategy would result in a global finite element space that is discontinuous between neighboring cells, but continuity can be restored by adding constraints that relate \glspl{dof} on neighboring cells.

\begin{figure}
  \def\Length{1}
  \def\Radius{0.03}
  \begin{tikzpicture}[scale=3]
    \LagrangeCell{0}{0}{\Length}{\Radius}{2}
      {{0,1,2,3,4,5,6,7,8}};
    \LagrangeCell{\Length}{0}{\Length}{\Radius}{4}
      {{9,10,11,12,13,14,15,16,17,18,19,20,21,22,23,24,25,26,27,28,29,30,31,32,33}};
  \end{tikzpicture}
  \hfill
  \begin{tikzpicture}[scale=3]
    \LagrangeCell{0}{0}{\Length}{\Radius}{2}
      {{0,1,2,3,4,5,6,7,8}};
    \LagrangeCell{\Length}{0}{\Length}{\Radius}{4}
      {{1,9,3,10,11,5,12,13,14,15,16,17,18,19,20,21,22,23,24,25,26,27,28,29,30}};
  \end{tikzpicture}%
  \caption{\it Enumeration of \gls{dof} indices on a mesh with two cells on which the left cell uses a $Q_2$ (bi-quadratic) Lagrange element and the right cell uses a $Q_4$ (bi-quartic) element. We distinguish between support points on vertices $(\bullet)$, lines $(\square)$ and quadrilaterals $(\circ)$.
  Left: Naive enumeration of \glspl{dof}. Continuity is ensured through constraints.
  Right: A better way in which we ``unify'' some \glspl{dof}.}
  \label{fig:naive-enumeration}
\end{figure}
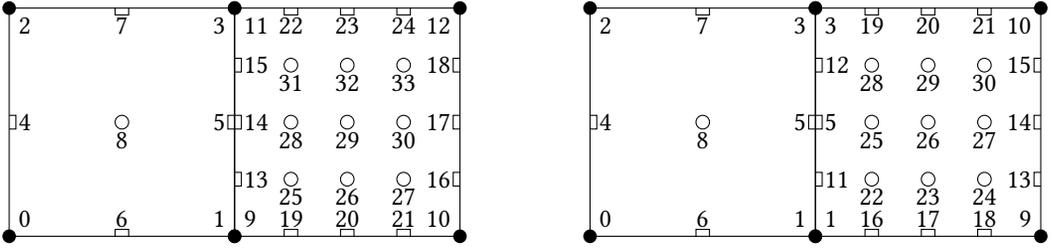

The left panel of Fig.~\ref{fig:naive-enumeration} illustrates this approach. Here, each cell's \glspl{dof} are independently enumerated. Continuity of the solution is then restored by introducing identity constraints of the form $U_9 = U_1$, $U_{11} = U_3$, $U_{14} = U_5$, in addition to the more traditional ``hanging node constraints'' $U_{13} = \tfrac 38 U_1 - \tfrac 18 U_3 + \tfrac 34 U_5$, $U_{15} = - \tfrac 18 U_1 + \tfrac 38 U_3 + \tfrac 34 U_5$.

While conceptually simple, this approach is wasteful as it introduces many more \glspl{dof} than necessary, along with a large number of constraints. In the extreme case of using $Q_1$ (tri-linear) Lagrange elements on all cells of a uniformly refined 3d mesh, one ends up with approximately eight times as many \glspl{dof}, $7/8$ of which are constrained. In actual test cases using \hp-adaptivity, \cite[Sec.~4.2]{bangerth2009} report that these ``unnecessary'' \glspl{dof} can be up to \SI{15}{\percent} of the total number of \glspl{dof} in 3d.

To avoid this wastefulness, the algorithms described in \cite{bangerth2009} ``unify'' \glspl{dof} where possible during the enumeration phase. For example, in the case shown in Fig.~\ref{fig:naive-enumeration}, the \glspl{dof} on shared vertices can be unified for the particular choice of elements adjacent to these vertices, as can be the \gls{dof} located on the common edge's midpoint. This leads to the enumeration shown on the right side of the figure, for which we then only need to add constraints 
$U_{11} = \tfrac 38 U_1 - \tfrac 18 U_3 + \tfrac 34 U_5$, $U_{12} = - \tfrac 18 U_1 + \tfrac 38 U_3 + \tfrac 34 U_5$.

At the same time, it is clear that this ``unification'' step requires knowing about the global indices of \glspl{dof} on neighboring cells \textit{during enumeration}, and this presents issues that need to be addressed in the parallel context if one of the cells adjacent to a vertex, edge, or face is a ghost cell. Furthermore, each process must know the active FE index not only for its locally owned cells, but also for ghost cells, before the enumeration can begin. We have to take into account all of these considerations in the extension of the algorithms of \cite{bangerth2009} to the parallel context in the next section.

\subsection{The parallel \hp-case}
\label{sec:enum-hp-parallel}

Having discussed the fundamental algorithms necessary to globally enumerate \glspl{dof} in the context of both parallel unstructured meshes, and for the sequential \hp-case, let us now turn to an algorithm that combines both of these features. As we will see, this algorithm turns out to be non-trivial.

\subsubsection{Goals for the parallel algorithm}

In developing such an enumeration algorithm, we are guided by the desire to come up with an enumeration that leads to a total number of \glspl{dof} that is independent of the number of processes. In other words, we do not want to treat vertices, edges, or faces that happen to lie on subdomain boundaries any different than if they were within the interior of a subdomain. We consider this an important feature to achieve scalable and predictable algorithms, and because it makes debugging problems easier. Furthermore, we would like to develop an algorithm that includes the ``unification step'' mentioned above to avoid generating too many trivial constraints.

At the end of the algorithm, each parallel process needs to know the globally unique indices of all \glspl{dof} located on the locally owned cells as well as on ghost cells, including the outer vertices, edges, and faces of ghost cells beyond which the current process has no knowledge of whether and how the mesh continues.

Finally, we want this algorithm to have linear complexity in the number of cells or the number of \glspl{dof}. We achieve this by stating it as a fixed-length series of loops over all cells owned by each process and, if necessary, over all ghost cells on this process. Because each process only loops over its own cells, and because the number of \glspl{dof} per process is balanced (see also Section~\ref{sec:load_balancing}), we obtain an algorithm that we expect to scale optimally both strongly and weakly.

The algorithm that achieves all of this -- see the discussion below -- can be broken down into seven distinct stages.
%whose work can be visually represented. For this visualization, we consider a simple example that consists of
In addition to their description, we illustrate each stage in an example for which we consider
a two-dimensional mesh of four neighboring cells meeting at a central vertex. On this mesh, we use bi-quadratic ($Q_2$) Lagrange elements with 9 unknowns on the bottom left and top right cell, and bi-quartic ($Q_4$) elements with 25 unknowns on the remaining two cells. Furthermore, we assume that the partitioning algorithm has divided the mesh into two subdomains: subdomain zero contains the bottom two cells, subdomain one the top two. This setup is shown in Fig.~\ref{fig:dof_enumeration}, where we illustrate the progress of the enumeration algorithm. Each figure shows the view from process zero on the left, and from process one on the right.

\subsubsection{Algorithm inputs}

The algorithm we describe in the following needs the following pieces of information as inputs:
\begin{itemize}
    \item A set of cells $K$ that constitute the \textit{locally
      owned} and \textit{ghost} cells, and information how neighboring
      cells are connected. The algorithm does not need to know where
      these cells are geometrically located in an ambient space --
      although this is of course important for the downstream
      application of the finite element method -- but only the
      topological connection of vertices, edges, and faces to the
      cells of which they are part of. The meshes we allow can be
      locally refined, with hanging nodes along edges. For practical
      reasons, we only allow a single hanging node per edge (i.e.,
      we enforce a 2:1 ratio of neighboring cell sizes), though this is
      not a critical restriction for the algorithm we will show.
    \item A process must know to which process each of its ghost cells belongs.
    %For each ghost cell on a process, that process needs to know which other process is the ``owner'' of that cell.
    Since we identify subdomain ids with process ranks in a \gls{mpi} universe, that means that we need to store the owning process's subdomain, or short owner, of all ghost cells.
    \item Each cell on every process has an associated \textit{global} identifier. This identifier is the same on all processes that store this cell, whether as part of their locally owned cells or as a ghost cell.
    \item For each locally owned or ghost cell, every process must know the active FE index --
    %For each cell, locally owned or ghost, each process needs to know the active FE index --
    that is, which element $\hat V_h^{a(K)}$ is in use on each cell $K$. Because the active FE index is typically computed only on each process's locally owned cells, this information needs to be exchanged between processes before the start of the algorithm; as with any ghost exchange of information, this can efficiently be done through point-to-point communication.
    \item For each element in the collection $\bigl\{\hat V_h^{(i)}\bigr\}_{i=1}^I$, the algorithm needs to know how many \glspl{dof} this element has per vertex, edge, face, or cell interior. For example, for the $Q_4$ Lagrange element in 2d that we use in our illustrative example, there is one \gls{dof} per vertex, three per edge, and nine per cell interior.
    \item A data structure that for each vertex, edge, face, and cell
      that is either locally owned or part of a ghost cell, can store
      the indices of \glspl{dof} located on this
      entity. Specifically, this implies that if one enumerates
      \glspl{dof} on one cell and sets their indices on that
      cell, that the indices of \glspl{dof} on vertices, edges,
      and faces are immediately also visible from neighboring cells
      without having to explicitly duplicate this information to other
      cells' index arrays. As discussed in \cite{bangerth2009},
      lower-dimensional entities (vertices, edges, and faces) have to
      be able to store multiple sets of indices, one for each of the
      finite element spaces used on the cells adjacent to the entity.
    \item For each pair of elements, the corresponding finite element
      implementations need to be able to identify whether two
      \glspl{dof} located on the same entity (vertex, edge, or face)
      can be unified. For example, our algorithm needs to be able to
      ask the combination of the $Q_2$ and $Q_4$ elements in the
      example shown in Fig.~\ref{fig:naive-enumeration} whether the
      single \gls{dof} each wants to store on a shared vertex can
      receive a single index, and whether the one $Q_2$ \gls{dof} on
      the shared edge can be unified with one of the three \glspl{dof}
      the $Q_4$ element wants to allocate on the common edge. In the
      example of Fig.~\ref{fig:naive-enumeration}, the answer is
      ``yes'' (\glspl{dof} 9, 11, and 14 in the left sub-figure can be
      identified with \glspl{dof} 1, 3, 5, and this is done in the
      right sub-figure). But this is not always the case, and
      consequently it is necessary to be able to query the adjacent finite elements about
      the possibility of \gls{dof} unification (for example,
      \glspl{dof} 13 and 15 in the left sub-figure cannot be unified with any
      \gls{dof} of the left cell, and similar cases also often happen
      when elements are defined with modal instead of nodal shape functions).
The details of how an answer to such a query can be implemented are not relevant to our description here, but are discussed at length in \cite{bangerth2009}.
\end{itemize}

\subsubsection{Description of the algorithm}

\let\oldthesubfigure\thesubfigure
\renewcommand{\thesubfigure}{After Stage \arabic{subfigure}}

\def\Length{1}
\def\Radius{0.03}

\begin{figure}
  \centering
  \begin{subfigure}{\textwidth}
    \begin{tikzpicture}[scale=3]
  \fill[color=Set1-F!80] (0, 0) rectangle (2*\Length, \Length);

  \tikzset{every node/.style={}};
  \LagrangeCell{0}{0}{\Length}{\Radius}{2}
    {{0,1,2,3,4,5,6,7,8}};
  \LagrangeCell{\Length}{0}{\Length}{\Radius}{4}
    {{9,10,11,12,13,14,15,16,17,18,19,20,21,22,23,24,25,26,27,28,29,30,31,32,33}};

  \tikzset{every node/.style={font=\ghostfont}};
  \LagrangeCell{0}{\Length}{\Length}{\Radius}{4} 
    {{"i","","i","i","i","i","i","i","i","i","i","i","i","i","i","i","i","i","i","i","i","i","i","i","i"}};
  \LagrangeCell{\Length}{\Length}{\Length}{\Radius}{2}
    {{"","i","i","i","i","i","i","i","i"}};
\end{tikzpicture}
    \hfill{}
    \begin{tikzpicture}[scale=3]
  \fill[color=Set1-F!80] (0, \Length) rectangle (2*\Length, 2*\Length);

  \tikzset{every node/.style={font=\ghostfont}};
  \LagrangeCell{0}{0}{\Length}{\Radius}{2}
    {{"i","i","i","","i","i","i","i","i"}};
  \LagrangeCell{\Length}{0}{\Length}{\Radius}{4}
    {{"i","i","","i","i","i","i","i","i","i","i","i","i","i","i","i","i","i","i","i","i","i","i","i","i"}};

  \tikzset{every node/.style={}};
  \LagrangeCell{0}{\Length}{\Length}{\Radius}{4}
    {{0,1,2,3,4,5,6,7,8,9,10,11,12,13,14,15,16,17,18,19,20,21,22,23,24}};
  \LagrangeCell{\Length}{\Length}{\Length}{\Radius}{2}
    {{25,26,27,28,29,30,31,32,33}};
\end{tikzpicture}
    \caption{Local enumeration.}
  \end{subfigure}
  \\[12pt]
  \begin{subfigure}{\textwidth}
   	\begin{tikzpicture}[scale=3]
  \tikzset{every node/.style={}};
  \LagrangeCell{0}{0}{\Length}{\Radius}{2}
    {{0,1,2,3,4,5,6,7,8}};
  \LagrangeCell{\Length}{0}{\Length}{\Radius}{4}
    {{9,10,11,12,13,14,15,16,17,18,19,20,21,22,23,24,25,26,27,28,29,30,31,32,33}};

  \tikzset{every node/.style={font=\ghostfont}};
  \LagrangeCell{0}{\Length}{\Length}{\Radius}{4}
    {{"i","","i","i","i","i","i","i","i","i","i","i","i","i","i","i","i","i","i","i","i","i","i","i","i"}};
  \LagrangeCell{\Length}{\Length}{\Length}{\Radius}{2}
    {{"","i","i","i","i","i","i","i","i","i","i","i","i","i","i","i"}};
\end{tikzpicture}
    \hfill{}
    \begin{tikzpicture}[scale=3]
  \fill[color=Set1-F!80] (\Length - 0.15, \Length) rectangle (\Length + 0.15, \Length + 0.15);

  \tikzset{every node/.style={font=\ghostfont}};
  \LagrangeCell{0}{0}{\Length}{\Radius}{2}
    {{"i","i","i","","i","i","i","i","i"}};
  \LagrangeCell{\Length}{0}{\Length}{\Radius}{4}
    {{"i","i","","i","i","i","i","i","i","i","i","i","i","i","i","i","i","i","i","i","i","i","i","i","i"}};

  \tikzset{every node/.style={}};
  \LagrangeCell{0}{\Length}{\Length}{\Radius}{4}
    {{0,"i",2,3,4,5,6,7,8,9,10,11,12,13,14,15,16,17,18,19,20,21,22,23,24}};
  \LagrangeCell{\Length}{\Length}{\Length}{\Radius}{2}
    {{"i",26,27,28,29,30,31,32,33}};
\end{tikzpicture}
    \caption{Tie-break.}
  \end{subfigure}
  
  \caption{\it Exemplary application of our enumeration algorithm for \glspl{dof}. Changes made at each step are highlighted. The left diagram of each subfigure depicts the situation for process 0, whereas the right side shows the domain from the perspective of process 1. The top half of each subfigure constitutes subdomain 1, while the bottom cells are assigned to subdomain 0. \Glspl{dof} on ghost cells are marked by italic indices.}
  \label{fig:dof_enumeration}
\end{figure}

\begin{figure}
  \ContinuedFloat
  \centering
  \begin{subfigure}{\textwidth}
    \begin{tikzpicture}[scale=3]
  \fill[color=Set1-F!80] (\Length - 0.15, 0) rectangle (\Length + 0.15, 0.15);
  \fill[color=Set1-F!80] (\Length - 0.15, 0.5*\Length - 0.1) rectangle (\Length + 0.15, 0.5*\Length + 0.1);
  \fill[color=Set1-F!80] (\Length - 0.15, \Length - 0.15) rectangle (\Length + 0.15, \Length);
  
  \fill[color=Set1-F!80] (1.5*\Length - 0.1, \Length - 0.13) rectangle (1.5*\Length + 0.1, \Length);
  \fill[color=Set1-F!80] (2*\Length - 0.15, \Length - 0.13) rectangle (2*\Length, \Length);

  \tikzset{every node/.style={}};
  \LagrangeCell{0}{0}{\Length}{\Radius}{2}
    {{0,1,2,3,4,5,6,7,8}};
  \LagrangeCell{\Length}{0}{\Length}{\Radius}{4}
    {{1,10,3,"i",13,5,15,16,17,18,19,20,21,22,"i",24,25,26,27,28,29,30,31,32,33}};

  \tikzset{every node/.style={font=\ghostfont}};
  \LagrangeCell{0}{\Length}{\Length}{\Radius}{4}
    {{"i","","i","i","i","i","i","i","i","i","i","i","i","i","i","i","i","i","i","i","i","i","i","i","i"}};
  \LagrangeCell{\Length}{\Length}{\Length}{\Radius}{2}
    {{"","i","i","i","i","i","i","i","i","i","i","i","i","i","i","i"}};
\end{tikzpicture}
    \hfill{}
    \begin{tikzpicture}[scale=3]
  \fill[color=Set1-F!80] (\Length - 0.15, 1.5*\Length - 0.1) rectangle (\Length + 0.15, 1.5*\Length + 0.1);
  \fill[color=Set1-F!80] (\Length - 0.15, 2*\Length - 0.15) rectangle (\Length + 0.15, 2*\Length);
  
  \fill[color=Set1-F!80] (0, \Length) rectangle (0.15, \Length + 0.13);
  \fill[color=Set1-F!80] (0.5*\Length - 0.1, \Length) rectangle (0.5*\Length + 0.1, \Length + 0.13);

  \tikzset{every node/.style={font=\ghostfont}};
  \LagrangeCell{0}{0}{\Length}{\Radius}{2}
    {{"i","i","i","","i","i","i","i","i"}};
  \LagrangeCell{\Length}{0}{\Length}{\Radius}{4}
    {{"i","i","","i","i","i","i","i","i","i","i","i","i","i","i","i","i","i","i","i","i","i","i","i","i"}};

  \tikzset{every node/.style={}};
  \LagrangeCell{0}{\Length}{\Length}{\Radius}{4}
    {{"i","i",2,27,4,5,6,7,29,9,10,"i",12,13,14,15,16,17,18,19,20,21,22,23,24}};
  \LagrangeCell{\Length}{\Length}{\Length}{\Radius}{2}
    {{"i",26,27,28,29,30,31,32,33}};
\end{tikzpicture}
    \caption{Unification.}
  \end{subfigure}
  \\[12pt]    
  \begin{subfigure}{\textwidth}
    \begin{tikzpicture}[scale=3]
  \fill[color=Set1-F!80] (0, 0) rectangle (2*\Length, \Length);
  
  \fill[white] (1.5*\Length - 0.1, \Length - 0.13) rectangle (1.5*\Length + 0.1, \Length + 0.13);
  \fill[white] (2*\Length - 0.15, \Length - 0.13) rectangle (2*\Length, \Length + 0.13);

  \tikzset{every node/.style={}};
  \LagrangeCell{0}{0}{\Length}{\Radius}{2}
    {{0,1,2,3,4,5,6,7,8}};
  \LagrangeCell{\Length}{0}{\Length}{\Radius}{4}
    {{1,9,3,"i",10,5,11,12,13,14,15,16,17,18,"i",19,20,21,22,23,24,25,26,27,28}};

  \tikzset{every node/.style={font=\ghostfont}};
  \LagrangeCell{0}{\Length}{\Length}{\Radius}{4}
    {{"i","","i","i","i","i","i","i","i","i","i","i","i","i","i","i","i","i","i","i","i","i","i","i","i"}};
  \LagrangeCell{\Length}{\Length}{\Length}{\Radius}{2}
    {{"","i","i","i","i","i","i","i","i","i","i","i","i","i","i","i"}};
\end{tikzpicture}
    \hfill{}
    \begin{tikzpicture}[scale=3]
  \fill[color=Set1-F!80] (0, \Length) rectangle (2*\Length, 2*\Length);
  
  \fill[white] (\Length - 0.15, \Length) rectangle (\Length + 0.15, \Length + 0.15);
  \fill[white] (0, \Length - 0.13) rectangle (0.15, \Length + 0.13);
  \fill[white] (0.5*\Length - 0.1, \Length - 0.13) rectangle (0.5*\Length + 0.1, \Length + 0.13);
  
  \tikzset{every node/.style={font=\ghostfont}};
  \LagrangeCell{0}{0}{\Length}{\Radius}{2}
    {{"i","i","i","","i","i","i","i","i"}};
  \LagrangeCell{\Length}{0}{\Length}{\Radius}{4}
    {{"i","i","","i","i","i","i","i","i","i","i","i","i","i","i","i","i","i","i","i","i","i","i","i","i"}};

  \tikzset{every node/.style={}};
  \LagrangeCell{0}{\Length}{\Length}{\Radius}{4}
    {{"i","i",29,50,30,31,32,33,52,34,35,"i",36,37,38,39,40,41,42,43,44,45,46,47,48}};
  \LagrangeCell{\Length}{\Length}{\Length}{\Radius}{2}
    {{"i",49,50,51,52,53,54,55,56}};
\end{tikzpicture}
	\caption{Global re-enumeration.}
  \end{subfigure}
  \caption{(continued) {\it Exemplary application of our enumeration algorithm for \glspl{dof}.}}
  % no label required
\end{figure}

\begin{figure}
  \ContinuedFloat
  \centering
  \begin{subfigure}{\textwidth}
    \begin{tikzpicture}[scale=3]
  \fill[color=Set1-F!80] (0,\Length) rectangle (2*\Length, 2*\Length);
  
  % \fill[white] (\Length - 0.15, \Length) rectangle (\Length + 0.15, \Length + 0.15);
  \fill[white] (0, \Length - 0.13) rectangle (0.15, \Length + 0.13);
  \fill[white] (0.5*\Length - 0.1, \Length - 0.13) rectangle (0.5*\Length + 0.1, \Length + 0.13);

  \tikzset{every node/.style={}};
  \LagrangeCell{0}{0}{\Length}{\Radius}{2}
    {{0,1,2,3,4,5,6,7,8}};
  \LagrangeCell{\Length}{0}{\Length}{\Radius}{4}
    {{1,9,3,"i",10,5,11,12,13,14,15,16,17,18,"i",19,20,21,22,23,24,25,26,27,28}};

  \tikzset{every node/.style={font=\ghostfont}};
  \LagrangeCell{0}{\Length}{\Length}{\Radius}{4}
    {{"i","",29,50,30,31,32,33,52,34,35,"i",36,37,38,39,40,41,42,43,44,45,46,47,48}};
  \LagrangeCell{\Length}{\Length}{\Length}{\Radius}{2}
    {{"",49,50,51,52,53,54,55,56}};
\end{tikzpicture}
    \hfill{}
    \begin{tikzpicture}[scale=3]
  \fill[color=Set1-F!80] (0,0) rectangle (2*\Length, \Length);
  
  \fill[white] (1.5*\Length - 0.1, \Length - 0.13) rectangle (1.5*\Length + 0.1, \Length + 0.13);
  \fill[white] (2*\Length - 0.15, \Length - 0.13) rectangle (2*\Length, \Length + 0.13);
  \fill[color=Set1-F!80] (\Length - 0.15, \Length) rectangle (\Length + 0.15, \Length + 0.15);

  \tikzset{every node/.style={font=\ghostfont}};
  \LagrangeCell{0}{0}{\Length}{\Radius}{2}
    {{0,1,2,"",4,5,6,7,8}};
  \LagrangeCell{\Length}{0}{\Length}{\Radius}{4}
    {{1,9,"","i",10,5,11,12,13,14,15,16,17,18,"i",19,20,21,22,23,24,25,26,27,28}};

  \tikzset{every node/.style={}};
  \LagrangeCell{0}{\Length}{\Length}{\Radius}{4}
    {{"i",3,29,50,30,31,32,33,52,34,35,"i",36,37,38,39,40,41,42,43,44,45,46,47,48}};
  \LagrangeCell{\Length}{\Length}{\Length}{\Radius}{2}
    {{3,49,50,51,52,53,54,55,56}};
\end{tikzpicture}
    \caption{Ghost exchange.}
  \end{subfigure}
      \\[12pt]
  \begin{subfigure}{\textwidth}
    \begin{tikzpicture}[scale=3]
  \fill[color=Set1-F!80] (0, \Length) rectangle (0.15, \Length + 0.13);
  \fill[color=Set1-F!80] (0.5*\Length - 0.1, \Length) rectangle (0.5*\Length + 0.1, \Length + 0.13);
  \fill[color=Set1-F!80] (1.5*\Length - 0.1, \Length - 0.13) rectangle (1.5*\Length + 0.1, \Length);
  \fill[color=Set1-F!80] (2*\Length - 0.15, \Length - 0.13) rectangle (2*\Length, \Length);
  % \fill[color=Set1-F!80] (\Length - 0.15, \Length - 0.15) rectangle (\Length + 0.15, \Length + 0.15);
  
  \tikzset{every node/.style={}};
  \LagrangeCell{0}{0}{\Length}{\Radius}{2}
    {{0,1,2,3,4,5,6,7,8}};
  \LagrangeCell{\Length}{0}{\Length}{\Radius}{4}
    {{1,9,3,49,10,5,11,12,13,14,15,16,17,18,54,19,20,21,22,23,24,25,26,27,28}};

  \tikzset{every node/.style={font=\ghostfont}};
  \LagrangeCell{0}{\Length}{\Length}{\Radius}{4}
    {{2,"",29,50,30,31,32,33,52,34,35,7,36,37,38,39,40,41,42,43,44,45,46,47,48}};
  \LagrangeCell{\Length}{\Length}{\Length}{\Radius}{2}
    {{"",49,50,51,52,53,54,55,56}};
\end{tikzpicture}
    \hfill{}
    \begin{tikzpicture}[scale=3]
  \fill[color=Set1-F!80] (0, \Length) rectangle (0.15, \Length + 0.13);
  \fill[color=Set1-F!80] (0.5*\Length - 0.1, \Length) rectangle (0.5*\Length + 0.1, \Length + 0.13);
  \fill[color=Set1-F!80] (1.5*\Length - 0.1, \Length - 0.13) rectangle (1.5*\Length + 0.1, \Length);
  \fill[color=Set1-F!80] (2*\Length - 0.15, \Length - 0.13) rectangle (2*\Length, \Length);
  % \fill[color=Set1-F!80] (\Length - 0.15, \Length - 0.15) rectangle (\Length + 0.15, \Length + 0.15);

  \tikzset{every node/.style={font=\ghostfont}};
  \LagrangeCell{0}{0}{\Length}{\Radius}{2}
    {{0,1,2,"",4,5,6,7,8}};
  \LagrangeCell{\Length}{0}{\Length}{\Radius}{4}
    {{1,9,"",49,10,5,11,12,13,14,15,16,17,18,54,19,20,21,22,23,24,25,26,27,28}};

  \tikzset{every node/.style={}};
  \LagrangeCell{0}{\Length}{\Length}{\Radius}{4}
    {{2,3,29,50,30,31,32,33,52,34,35,7,36,37,38,39,40,41,42,43,44,45,46,47,48}};
  \LagrangeCell{\Length}{\Length}{\Length}{\Radius}{2}
    {{3,49,50,51,52,53,54,55,56}};
\end{tikzpicture}
    \caption{Merge on interfaces.}
  \end{subfigure}
  \caption{(continued) {\it Exemplary application of our enumeration algorithm for \glspl{dof}.}}
  % no label required
\end{figure}

\renewcommand{\thesubfigure}{\oldthesubfigure}

The algorithm we show consists of seven stages (plus an initial memory allocation and initialization stage), as detailed below. To make understanding it easier, we illustrate each step in Fig.~\ref{fig:dof_enumeration} for our model test case; note that the figure is continued over several pages. In our description, we follow the same nomenclature as \cite{bangerth2012}. Specifically, we generally use an index $p$ or $q$ for subdomains (identified with process ranks in a \gls{mpi} universe), and we denote the set of all \textit{locally owned cells} on process $p$ by $\cellsp{loc}$, the set of all \textit{ghost cells} by $\cellsp{ghost}$, and the set of all \textit{locally relevant cells} by $\cellsp{rel} = \cellsp{loc} \cup \cellsp{ghost}$. In the description below, the process index $p$ is generically used to identify the ``current'' process, i.e., the rank of the process that is executing the algorithm.

Our algorithm then proceeds in the following steps:

\begin{enumerate}
\setcounter{enumi}{-1}
  \item \textbf{Initialization (without illustration).}
  % no communication
  Loop over all locally relevant cells $K \in \cellsp{rel}$, and on each of its vertices, edges, faces, and $K$ itself allocate enough space to store as many \gls{dof} indices as are necessary for the element identified by the active FE index $a(K)$. If a neighboring element has already allocated space for the same active FE index, then no additional space is necessary. In other words, for each entity within $\cellsp{rel}$, we need to allocate space for a map from the active FE indices of adjacent cells to an array of indices of \glspl{dof} indices.
  
  Once space is allocated, all \gls{dof} indices are set to an invalid value that we denote by $i$ in the following (for example $i \coloneqq -1$).

  \textit{Communication:} This stage does not require any communication.

  \item \textbf{Partition-local enumeration.}
  % no communication
  Iterate over all locally owned cells $K \in \cellsp{loc}$. For each of the vertices, edges, faces, and the cell interior, assign valid \gls{dof} indices in ascending order, starting from zero, if indices have not already been assigned for an entity and the current $a(K)$.

  \textit{Communication:} This stage does not require any communication.
  
  \item \textbf{Tie-break.}
  % no communication
  Iterate over all locally owned cells $K \in \cellsp{loc}$. If a vertex, edge, or face that is part of $K$ is also part of an adjacent ghost cell $K'_\text{ghost}$ so that $a(K)=a(K'_\text{ghost})$, and if $K'_\text{ghost}$ belongs to a subdomain of lower rank $q < p$, then invalidate all \glspl{dof} on this mesh entity by setting their index to the invalid value $i$.

  \textit{Communication:} This stage does not require any communication.
  
  \item \textbf{Unification.}
  % no communication
  Iterate over all locally owned cells $K \in \cellsp{loc}$. For all shared \glspl{dof} on vertices, edges, and faces to neighboring cells $K'$ (locally owned or ghost), ask the elements corresponding to active FE indices $a(K)$ and $a(K')$ whether some of the \glspl{dof} can be unified between the two elements.
  If $K'$ is also a locally owned cell, perform the unification by replacing one index (or a set of indices) by the corresponding index of the other \gls{dof} to which it is unified.
  If $K'$ is a ghost cell, and if the \gls{dof} on $K$ needs to be unified with the corresponding one on $K'$ (rather than the other way around), then set the index of the \gls{dof} on $K$ to the invalid value $i$.

  \textit{Communication:} This stage does not require any communication.
\end{enumerate}

At this point in the algorithm, each process knows which \glspl{dof} are owned by this process -- namely, the ones on locally owned cells that are enumerated as anything other than $i$ -- although the final indices of these \glspl{dof} are not yet known.

\begin{enumerate}[resume]
  \item \textbf{Global re-enumeration.}
  % one MPI_(Ex)scan
  Iterate over all locally owned cells $K \in \cellsp{loc}$ and
  re-enumerate those \gls{dof} indices in ascending order that have a
  valid value assigned, ignoring all invalid indices. Store the total
  number of all valid \gls{dof} indices on this subdomain as $n_p$. In a next step, shift all indices by the number of \glspl{dof} that are owned by all processes of lower rank $q < p$, or in other words, by $\sum_{q=0}^{p-1} n_q$. Computing this shift corresponds to a prefix sum or exclusive scan, and can be obtained via \texttt{MPI\_Exscan} \cite{mpi40}.

  \textit{Communication:} This stage requires one global \texttt{MPI\_Exscan} operation on a single integer of sufficient size to hold the largest \gls{dof} index.
\end{enumerate}

At this stage, each process has (consecutively) enumerated a certain subset of \glspl{dof}, and we call these the ``locally owned \glspl{dof}''. In later use, each process then owns the corresponding rows of matrices and entries in vectors, but the concept of locally owned \glspl{dof} is otherwise of no importance to the remainder of the algorithm. Importantly, however, we still need to ensure that each process learns of the remaining \glspl{dof} that are located on locally owned or ghost cells and whose indices are not currently known.

\begin{enumerate}[resume]
  \item \textbf{Ghost exchange.}
    % communication as described
    In this step, we need to send sets of
    indices from those locally owned cells $K \in \cellsp{loc}$ that
    are ghost cells on other processes, to those processes on which they are ghost cells. We do this in the following steps:
  \begin{enumerate}[label=\alph*.]
    \item For each process $q\neq p$ that is adjacent to $p$, allocate
      a map with keys corresponding to global cell identifiers and
      values equal to a list of indices of those \glspl{dof} defined on this cell.
    \item Iterate over all $K \in \cellsp{loc}$. If $K$ is a ghost cell on process $q$, then add the global identifier of $K$ and the list of \glspl{dof} on $K$ to the map for process $q$.
    \item Send all of the maps to their designed process $q$ via nonblocking point-to-point communication (e.g., using \texttt{MPI\_Isend} \cite{mpi40}).
    \item Receive data containers from processes of adjacent subdomains $q$ via nonblocking point-to-point communication (e.g., using \texttt{MPI\_Irecv} \cite{mpi40}). The data so received corresponds to the \gls{dof} indices on all ghost cells of this subdomain $p$. On each of these cells, set the received \gls{dof} indices accordingly.
  \end{enumerate}
  All communication in this step is symmetric, which means that a process only receives data from another process when it also sends data to it. Thus, there is no need to negotiate communication.

  \textit{Communication:} This stage requires point-to-point communication between all processes that own neighboring sub-domains. The data sent consists of the indices of \glspl{dof} on all those cells that are locally owned by the sending process and that are ghost cells on the receiving process.
\end{enumerate}
After this ghost exchange, each \gls{dof} on an interface between a locally owned and a ghost cell has exactly one valid index assigned.
\begin{enumerate}[resume]
  \item \textbf{Merge on interfaces.}
  % no communication
  Iterate over all locally relevant cells $K \in \cellsp{rel}$. On interfaces between locally owned and ghost cells, set all remaining invalid \gls{dof} indices to the corresponding valid one.

  \textit{Communication:} This stage does not require any communication.
\end{enumerate}
At this stage, all processes know the correct indices for all \glspl{dof} located on locally owned cells. However, during the ghost exchange in stage 5 above, some processes may have sent index sets for some cells that may still contain the invalid index $i$ and not all of these can be resolved through unification with locally known indices in stage 6. This is not illustrated in the figures but would require a larger example mesh; the source of these $i$ markers are if a ghost cell owned by process $q$ does not only border a cell owned by process $p$, but also a cell owned by yet another process $q'$ that is not a neighbor of $p$, and process $q$ will only learn about indices on this cell \textit{as part of the ghost exchange with $q'$ itself}. As a consequence, we have to repeat stage 5 one more time:
\begin{enumerate}[resume]
  \item \textbf{Ghost exchange (without illustration).}
  % communication as described in stage 5
  Repeat the steps of stage 5. However this time, only data from those cells have to be communicated which had invalid \gls{dof} indices prior to stage 5d.

  \textit{Communication:} This stage requires point-to-point communication between all processes that own neighboring sub-domains. The data sent consists of the indices of \glspl{dof} on a subset of all those cells that are locally owned by the sending process and that are ghost cells on the receiving process.
\end{enumerate}

At the end of this algorithm, all global \gls{dof} indices have been
set correctly, and every process knows the indices of \glspl{dof}
located on locally owned and ghost cells (i.e., the ``locally relevant
\glspl{dof}'' in the terminology of \cite{bangerth2012}). The proof that this is so
is given by the following three considerations: (i)~At the end of
stage 4, we know the final indices of all \glspl{dof} on the
locally owned cells with the exception of ones on interfaces to ghost
cells if these \glspl{dof} are owned by another process. (ii)~In stage 5 of
the algorithm, through the first ghost exchange, every process
receives information about \gls{dof} indices from the owners of these
DoFs. Stage 6, which merges \glspl{dof} on interfaces, then completes the
knowledge of all \glspl{dof} on locally owned cells. (iii)~Because now every
process knows about all indices for the \glspl{dof} on locally owned cells,
in stage 7, during the second ghost exchange, every process receives
information from the owners of all ghost cells that allows completing
all \gls{dof} indices on ghost cells.

\begin{remark}
In three-dimensional scenarios, \cite[Sec.~4.6]{bangerth2009} points out possible complications with circular constraints during \gls{dof} unification whenever three or more different finite elements share a common edge. We have not found other satisfactory solutions for this problem in the intervening years, and consequently continue to implement the suggestion in \cite{bangerth2009}: all \glspl{dof} on such edges are excluded from the unification step and will be treated separately via constraints. Since the decision to use or not use the unification algorithm on these edges is independent of whether the adjacent cells are on the same or different processes, this decision has no bearing on our overall goals of ensuring that the number of used indices be independent of the partitioning of the mesh. In the examples presented in Section~\ref{sec:results}, the fraction of identity constraints stays below \SI{3}{\percent}.
\end{remark}

\begin{remark}
During stage 3 of our example in Fig.~\ref{fig:dof_enumeration}, we follow the \gls{dof} unification procedure as described in \cite{bangerth2009}: if different finite elements meet on a subdomain interface, all shared \glspl{dof} will be assigned to the finite element representing the common function space (that is, when using elements within the same family, the one with the lower polynomial degree). Of course, different decisions are possible, which might have an impact on parallel performance. For example, \cite[Remark 2]{bangerth2012} pointed out that on a face, all \glspl{dof} should belong to the same subdomain to speed up parallel matrix-vector multiplications. We implemented such an enumeration algorithm as an alternative to the one presented here. For the Laplace example used for the weighted load balancing experiments described in Section~\ref{sec:results-load-balancing}, we found that both implementations take the same run time (\textless{} \SI{1}{\percent} deviation).
\end{remark}

% ----------------------
% --- Load balancing ---
% ----------------------

\section{Load balancing}
\label{sec:load_balancing}

In order to enable our algorithms to scale well, we need to ensure
that each process does roughly the same amount of work. In contrast to
\h-adaptively refined meshes, a major difficulty here is that the
workload per cell is not the same: different parts of the overall \hp-adaptive algorithm scale differently with the number $n_\text{\glsfmtshortpl{dof}}$ of unknowns per cell -- for example, the cost of enumerating \glspl{dof} on a cell is proportional to $n_\text{\glsfmtshortpl{dof}}$, whereas assembling cell-local contributions to the global system costs ${\mathcal O}(n_\text{\glsfmtshortpl{dof}}^3)$, unless one uses specific features of the finite element basis functions. More importantly, how the cost of a linear solver or algebraic multigrid implementation -- together the largest contribution to a program's run time -- scales with the polynomial degree or number of unknowns on a cell is quite difficult to estimate \textit{a priori}. As a consequence, when using different polynomial degrees on different cells, it is not easy to derive theoretically what the computational \textit{cost} of a cell is going to be, and consequently how to weigh each cell.

\cite{oden1994, patra1995} investigate different decomposition and
load balancing strategies with various types of weights, which are
closely tied to their \hp-adaptive algorithm. These studies use the
number of \glspl{dof} on a given cell as a natural choice for the weight of that cell, but we believe that this does not reflect the computational effort accurately for the reasons pointed out above.
Herein, we instead use an empirical approach in which we assume that the relative cost $w$ of a cell $K$ can be expressed as
\begin{align}
\label{eq:load-balance-weight}
    w(K) = n_\text{\glsfmtshortpl{dof}}(K)^c,
\end{align}
with some, \textit{a priori} unknown, exponent $c$. During load balancing, we then weigh each cell with this factor and seek to partition meshes so that the sums of weights of the cells in each partition are roughly equal. For this purpose, we use the partitioning algorithms provided by the \pforest{} library that are described in \cite[Sec.~3.3]{burstedde2011}.

We experimentally determine the value for the exponent $c$ for
which the overall run time of our program is minimized, and show
results to this end in Section~\ref{sec:results-load-balancing}. From
the considerations above, one would expect that the minimum should be in the range $1\le c\le 3$, and this indeed turns out to be the case.

It is worth mentioning that the approach only minimizes the
\textit{overall} run time, but likely leaves each individual operation
sub-optimally load balanced. This imbalance is a common problem when a
program executes algorithms whose cell-local costs are not
proportional (see, for example, \cite{gassmoller2018}) and can only be
solved by re-partitioning data structures between the different phases
of a program -- say, between matrix assembly and the actual solver
phase. Exploring this issue is beyond the scope of our study;
moreover, as we will show in Section~\ref{sec:results}, for the test
cases we consider, the solver's contribution to the overall run time
is so dominant that it is not worth trying to better load-balance any
of the non-solver components.

% ----------------------------
% --- Global data transfer ---
% ----------------------------

\section{Packing, unpacking, and transferring data}
\label{sec:data_transfer}

A frequent operation in finite element codes is the serialization of all information associated with a cell into an array, and moving this data. Examples for where this operation is relevant are re-partitioning a mesh among processes after refinement, and the generation of checkpoints for later restart. In such cases, it is often convenient to write all information associated with the cells of one process into contiguous buffers.

For \h-adaptive meshes, this presents few challenges since the size of
the data associated with every cell is the same and, consequently, can be packed into buffers of fixed size per cell. On the other hand, for \hp-methods, different cells require different buffer sizes for efficiency, and creating contiguous storage schemes for all data on each process requires a bit more thought. Thought is also necessary when devising mechanisms to subsequently transfer this data to other processes.

In practice, we implement such schemes using a two-stage process: in a first stage, we assess how much memory the data on each cell requires, and allocate a contiguous array that can hold information from all cells. In this phase, we also build a second array that holds the offsets into the first array at which the data from each cell starts. The storage scheme therefore resembles the way sparse matrices are commonly stored in \gls{csr}. In a second stage, we copy the actual data from each cell into the respective part of the array.

For serialization, one can then write the two arrays in their entirety
to disk. For re-partitioning, parts of the arrays have to be sent to
different processes based on which process will own a cell. For this
step, it is useful to sort the order in which cells are represented in
the two arrays in such a way that data destined for one target is
stored as one contiguous part of the arrays. In this way, all
information to be sent to one process can be transferred with a single
non-blocking point-to-point send operation (with non-uniform buffer
sizes) for each of the two arrays,
without the need for further copy operations. How this can be
efficiently implemented is described in detail in \cite[Sec.~5.2]{burstedde2020}.

%\begin{figure}
%  \centering
%  \input{figures/transfer/memory_for_transfer}%
%  \todo[inline]{The figure is pretty, but I'm not sure it provides sufficient additional information not represented in the text. If we want to keep it, we ought to make sure that the data stored is really of variable size, and we should have a depiction of the second array that stores offsets.}
%  \caption{\it Division of contiguous chunk of memory for data transfer.}
%  \label{fig:memory}
%\end{figure}

% -------------------------
% --- Numerical results ---
% -------------------------

\section{Numerical Results}
\label{sec:results}

Ultimately, the algorithms we have presented in Sections~\ref{sec:enumeration}, \ref{sec:load_balancing}, and \ref{sec:data_transfer} are only useful if they can be efficiently implemented. In this section, we assess our approaches using two test cases: a two-dimensional Laplace problem, and a three-dimensional Stokes problem. We discuss these in Section~\ref{sec:testcases} below.

Based on these test problems, we first assess how one needs to choose load balancing weights for each cell based on the polynomial degree of the finite element applied (Section~\ref{sec:results-load-balancing}). Using the resulting load balancing strategy, we then discuss how our algorithms scale in Section~\ref{sec:results-scaling}; an important question to discuss in this context will be how one would actually define and measure ``scalability'' in the context of \hp-adaptive methods.

All of the results shown in this section have been obtained using codes that are variations of tutorial programs of the \dealii{} library. All features discussed in this paper are implemented in \dealii{}, see also \cite{dealii94,arndt2021}. All data were generated with the tool \hpbox{} \cite{fehling2022}.

\subsection{Test cases}
\label{sec:testcases}

We evaluate the performance of our algorithms using two test cases discussed below: a two-dimensional Laplace equation posed on the L-shaped domain, and a three-dimensional Stokes problem posed on a domain that resembles a forked (``Y''-shaped) pipe. Both of these cases are chosen because the domain induces corner singularities in the solution, resulting in parts of the domain where either large cells with high-order elements or small cells with low-order elements are best suited to approximate the exact solution. In other words, these cases mimic practical situations that are well suited to \hp-adaptive methods. Furthermore, being able to demonstrate our algorithms on both a relatively simple, scalar two-dimensional problem and a much more complex three-dimensional, coupled vector-valued problem illustrates the range and limitations of our algorithms.

In each test case, we start from a coarse discretisation of the problem, solve it, and refine it in multiple iterations to end up with a mesh tailored to the problem. For this purpose, we need mechanisms to decide which cells we want to refine and how. We use an error estimator based on \cite{kelly1983} to mark cells for general refinement. Further, we use a smoothness estimator based on the decay of Legendre coefficients as described by \cite{eibner2007, houston2005, mavriplis1994} to decide how we want to refine each cell. We employ fixed number refinement for both \h- and \p-refinement, which means that the fraction of cells we are going to refine is always the same. We state our choice of fractions in the descriptions below. The mesh is repartitioned after each refinement iteration.

\subsubsection{Test case 1: A Laplace problem on the L-shaped domain}
\label{sec:test-case-1}

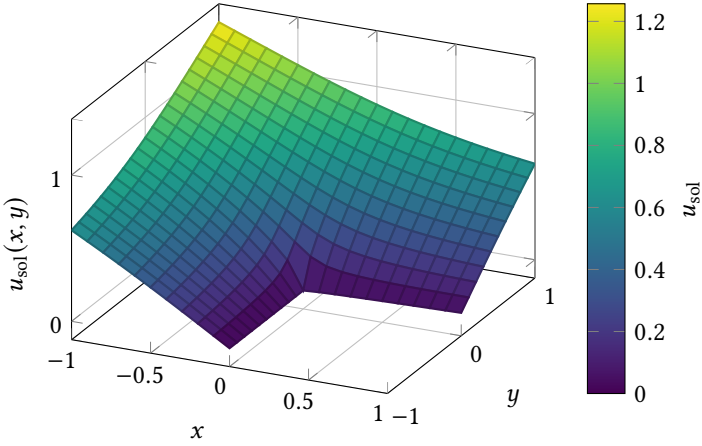
\begin{figure}
  \centering
  \begin{tikzpicture}
\begin{axis}[
  scale=0.9,
  grid=major,
  colorbar,
  colormap/viridis,
  colorbar style={
    ylabel={$u_\mathrm{sol}$},
  },
  xlabel=$x$,
  ylabel=$y$,
  zlabel={$u_\mathrm{sol}(x,y)$}
  ]
\addplot3[
  surf,
  domain=-1:1,
  y domain=-1:1,
  restrict expr to domain={(x>0)&&(y<0)}{0:0},
  samples=21 %11
  ]{pow(x*x+y*y,0.33)*sin(0.66*(atan2(y,-x)+90)))};
\end{axis}
\end{tikzpicture}
  \caption{\it The solution \eqref{eq:ficherasolution} of the Laplace problem \eqref{eq:laplaceproblem}  on the L-shaped domain.}
  \label{fig:laplace_solution}
\end{figure}

Our first test case concerns the solution of the Laplace problem with Dirichlet boundary conditions:
\begin{align}
  \label{eq:laplaceproblem}
  - \Delta u(\vec{x}) &= 0 \quad\text{on}\quad \Omega \,\text{,} & u(\vec{x}) &= u_\text{sol}(\vec{x}) \quad\text{on}\quad \partial\Omega,
\end{align}
where we choose $\Omega\subset\mathbb{R}^2$ as the L-shaped domain, $\Omega=(-1,1)^2\backslash[0,1]\times[-1,0]$.
It is well understood that on such domains, the Laplace equation admits a singular solution; indeed, in polar coordinates $r = \sqrt{x^2 + y^2} > 0$ and $\theta = \arctan (y/x)$, the function
\begin{align}
  \label{eq:ficherasolution}
  u_\text{sol}(\vec{x}) &= r^\alpha \sin(\alpha \, \theta)
\end{align}
is a solution for an opening angle at the reentrant corner of $\pi / \alpha$ with $\alpha \in (1/2,1)$. For the L-shaped domain, we have $\alpha = 2/3$ and the corresponding solution is shown in Fig.~\ref{fig:laplace_solution}. We impose $u=u_\text{sol}$ as the boundary condition on $\partial\Omega$, and the resulting (exact) solution of the Laplace equation that we seek to compute is then $u=u_\text{sol}$ everywhere in $\Omega$.

This solution is singular at the origin: with unit vectors $\vec{e}_r = \cos(\theta) \vec{e}_x - \sin(\theta) \vec{e}_y$ and $\vec{e}_\theta = \sin(\theta) \vec{e}_x + \cos(\theta) \vec{e}_y$, we find that
\begin{align}
  \nabla u_\text{sol}(\vec{x}) &= \alpha r^{\alpha - 1} \left[ \sin(\alpha \, \theta) \vec{e}_r + \cos(\alpha \, \theta) \vec{e}_\theta \right],
\end{align}
and consequently $\lim\limits_{r \rightarrow 0} \left\| \nabla
u_\text{sol}(\vec{x}) \right\|_{2} = \infty$ for our choice of $\alpha$.

The numerical solution of the Laplace equation on the L-shaped domain is a classical test case. For example, \cite{mitchell2014} presents several benchmarks for \hp-adaptation for this situation. 
A similar scenario is also used in the step-75 tutorial of the \dealii{} library \cite{fehling2021}.

In our study, we choose Lagrange elements $Q_k$ with polynomial degrees $k = 2,\dots,7$. We
% create the mesh in Fig.~\ref{fig:laplace_approximation} by marking
mark \SI{30}{\percent} of cells for refinement and \SI{3}{\percent} for coarsening, from which we pick \SI{90}{\percent} to be \p-adapted and \SI{10}{\percent} to be \h-adapted. We choose to favor \p-refinement since the only non-smooth part of the solution is around the point singularity at the origin.

\begin{figure}
  \begin{tikzpicture}
\begin{axis}[
  scale=0.92,
  xmin=-1,xmax=1,
  ymin=-1,ymax=1,
  unit vector ratio={1 1},
  tick align=outside,
  xtick={-1,-0.5,0,0.5,1},
  ytick={-1,-0.5,0,0.5,1},
  xlabel=$x$,
  ylabel=$y$,
  colormap/OrRd,
  colorbar sampled,
  colorbar horizontal,
  colorbar style={xlabel={polynomial degree $k$}, xtick={2,3,...,7}, samples=7},
  point meta min=1.5,
  point meta max=7.5
]

\addplot graphics [
  xmin=-1,xmax=1,
  ymin=-1,ymax=1,
] {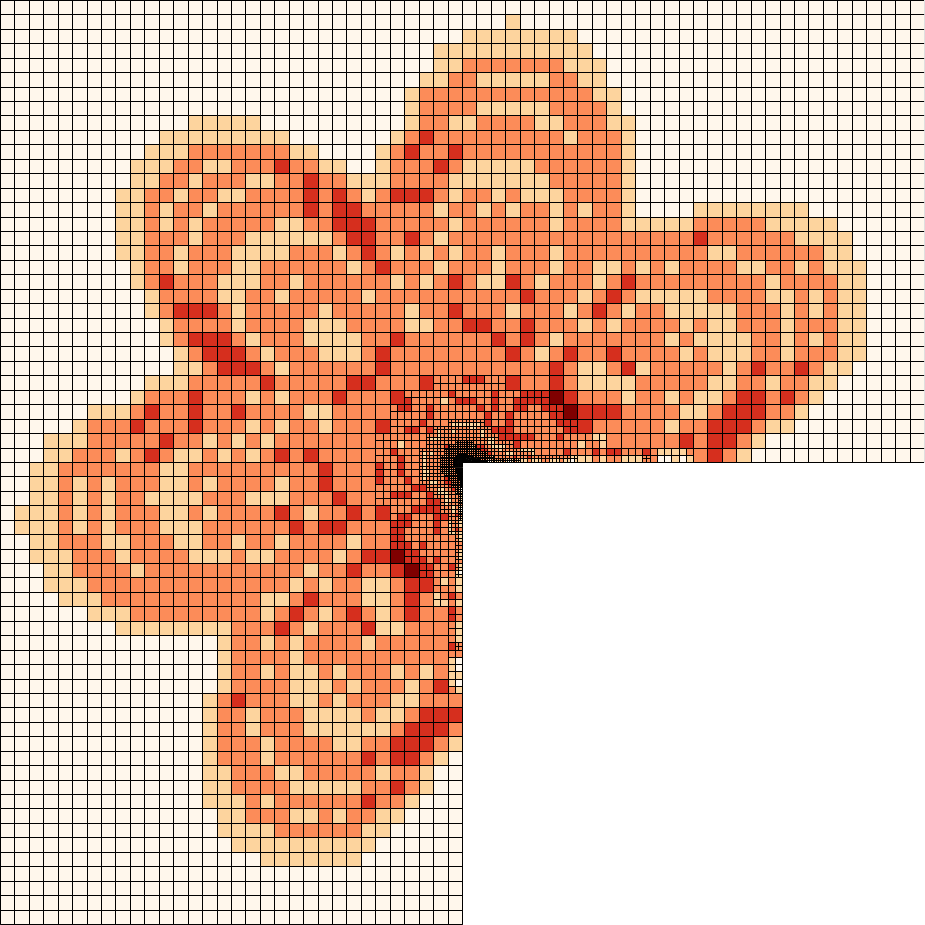};
\end{axis}
\end{tikzpicture}
  \hfill
  \begin{tikzpicture}
\begin{axis}[
  scale=0.92,
  xmin=-1,xmax=1,
  ymin=-1,ymax=1,
  unit vector ratio={1 1},
  tick align=outside,
  xtick={-1,-0.5,0,0.5,1},
  ytick={-1,-0.5,0,0.5,1},
  xlabel=$x$,
%   ylabel=$y$,
  colormap/Paired-12,
  colorbar sampled,
  colorbar horizontal,
  colorbar style={xlabel={subdomain id}, xtick={0,1,...,11}, samples=13},
  colormap access=piecewise const,
  point meta min=-0.5,
  point meta max=11.5
]

\addplot graphics [
  xmin=-1,xmax=1,
  ymin=-1,ymax=1,
] {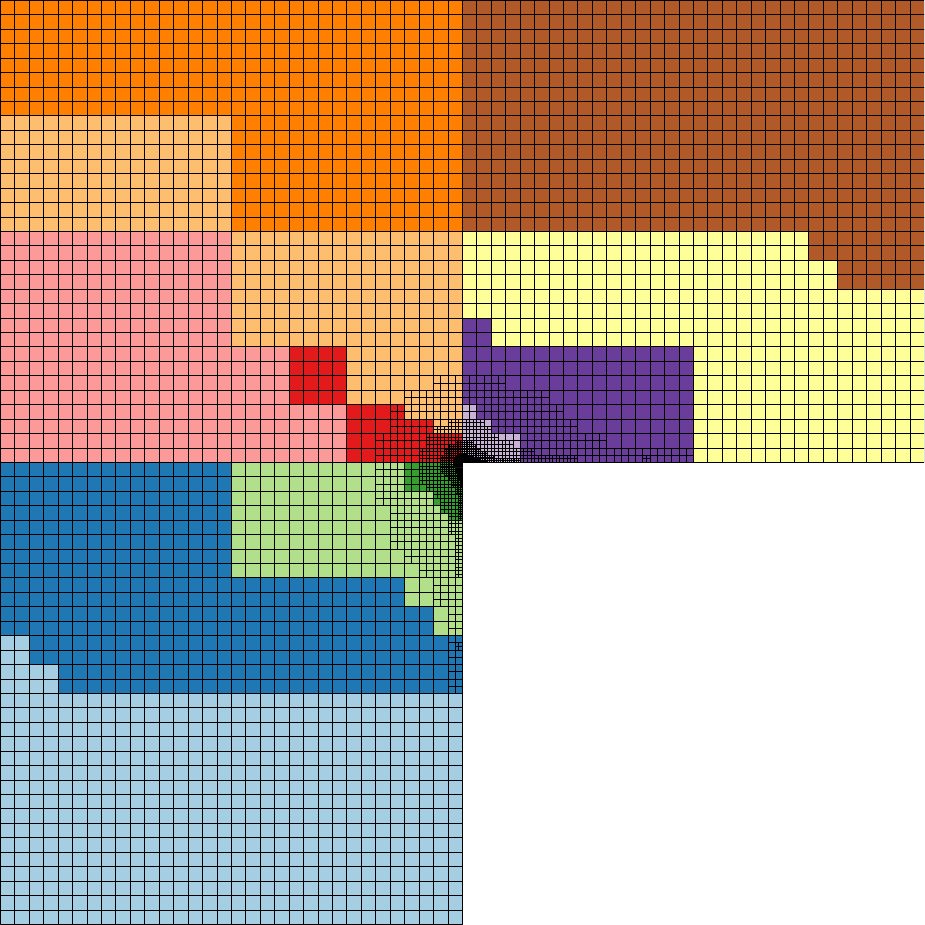};
\end{axis}
\end{tikzpicture}
  \caption{\it Numerical approximation of the Laplace problem \eqref{eq:laplaceproblem} after six adaptation cycles and five initial global refinements. Left: The mesh and polynomial degrees used on each cell. Right: Partitioning of the mesh onto 12 \gls{mpi} processes with a load balancing weighting exponent of $c = 1.9$.}
  \label{fig:laplace_approximation}
\end{figure}

Fig.~\ref{fig:laplace_approximation} shows a typical \hp-mesh and its partitioning from a sequence of adaptive refinements.
It illustrates that the corner singularity requires \h-adaptation resulting in small cells, whereas further away from the corner, the solution is smoother and can be resolved on relatively coarse meshes using high polynomial degrees. Far away from the origin, the estimated errors are low so that large cells and low polynomial degrees are sufficient.
The lobe pattern results from the anisotropic resolution property of polynomials on quadrilaterals.

The numerical scheme we choose to solve this problem is based on \trilinos{} \cite{heroux2005} for parallel linear algebra, and uses the \ml{} package \cite{gee2007} as an \gls{amg} preconditioner inside a conjugate gradient iteration.

\subsubsection{Test case 2: Flow through a Y-pipe}
\label{sec:test-case-2}

\begin{figure}
  \begin{minipage}{0.49\textwidth}
    \includegraphics[width=\textwidth]{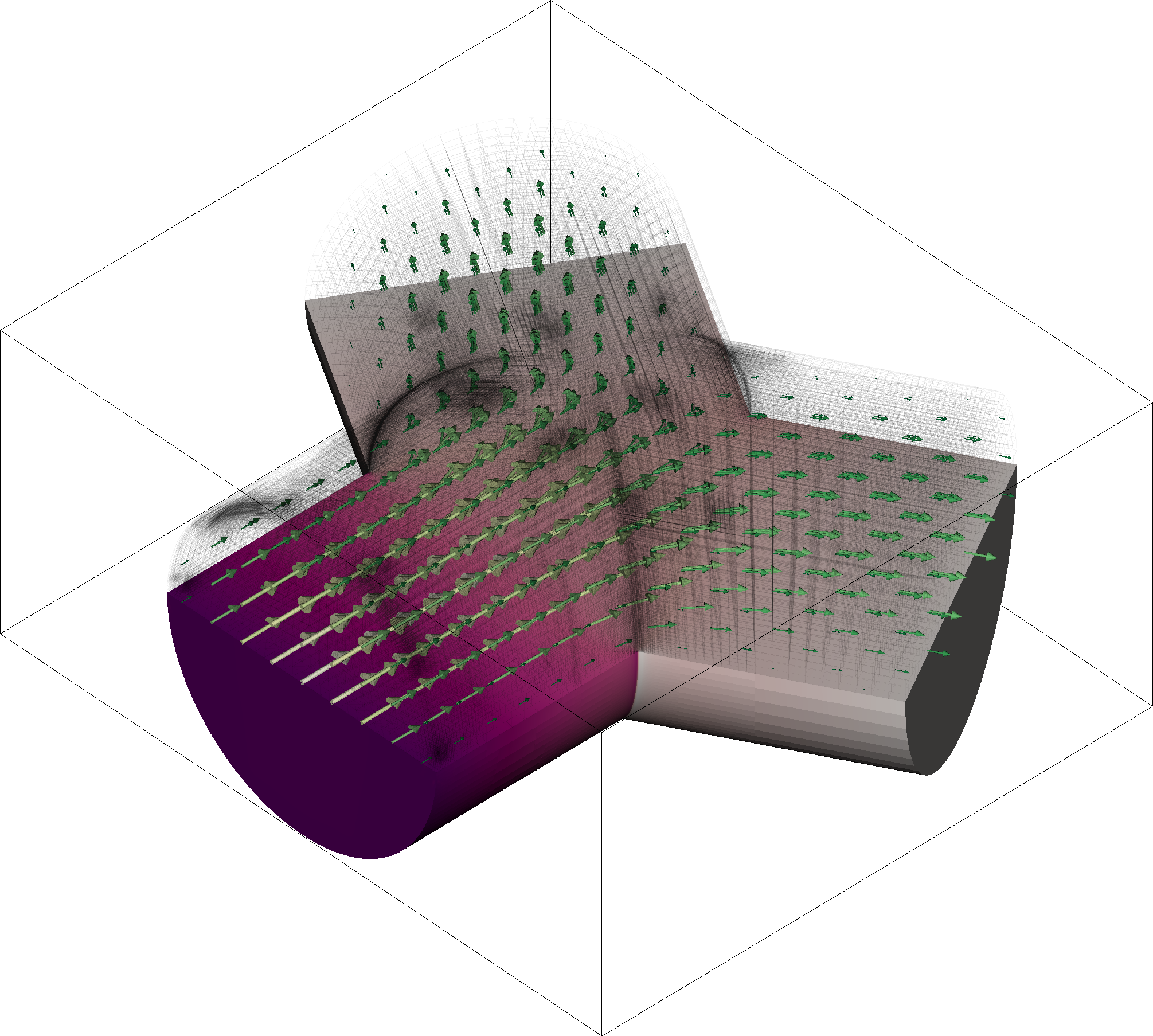}
  \end{minipage}
  \hfill
  \begin{minipage}{0.49\textwidth}
    \hfill
\adjustbox{valign=c}{
  \begin{tikzpicture}
    \pgfplotscolorbardrawstandalone[ 
      colormap/RdPu,
      point meta min=0,
      point meta max=1.233,
      colorbar style={
        ylabel={pressure},
        height=5cm,
        ytick={0,0.2,...,1,1.233},
      }
    ]
  \end{tikzpicture}
}
\hfill
\adjustbox{valign=c}{
  \begin{tikzpicture}
    \pgfplotscolorbardrawstandalone[ 
      colormap/YlGn,
      colormap={GnYl}{
        indices of colormap={
          \pgfplotscolormaplastindexof{YlGn},...,0 of YlGn}
      },
      point meta min=0,
      point meta max=1,
      colorbar style={
        ylabel={velocity},
        height=5cm,
      }
    ]
  \end{tikzpicture}
}
\hfill
\adjustbox{valign=c}{
  \begin{tikzpicture}
    \pgfplotscolorbardrawstandalone[
      colormap/RdYlBu,
      colormap={BuYlRd}{
        indices of colormap={
          \pgfplotscolormaplastindexof{RdYlBu},...,0 of RdYlBu}
      },
      colormap access=piecewise const,
      point meta min=2.5,
      point meta max=6.5,
      colorbar sampled,
      colorbar style={
        ylabel={polynomial degree $k$ of $\pmb{Q}_k/Q_{k-1}$},
        height=5cm,
        ytick={3,4,5,6},
        samples=5,
      }
    ]
  \end{tikzpicture}
}
  \end{minipage}\\
  \begin{minipage}{0.49\textwidth}
    \includegraphics[width=\textwidth]{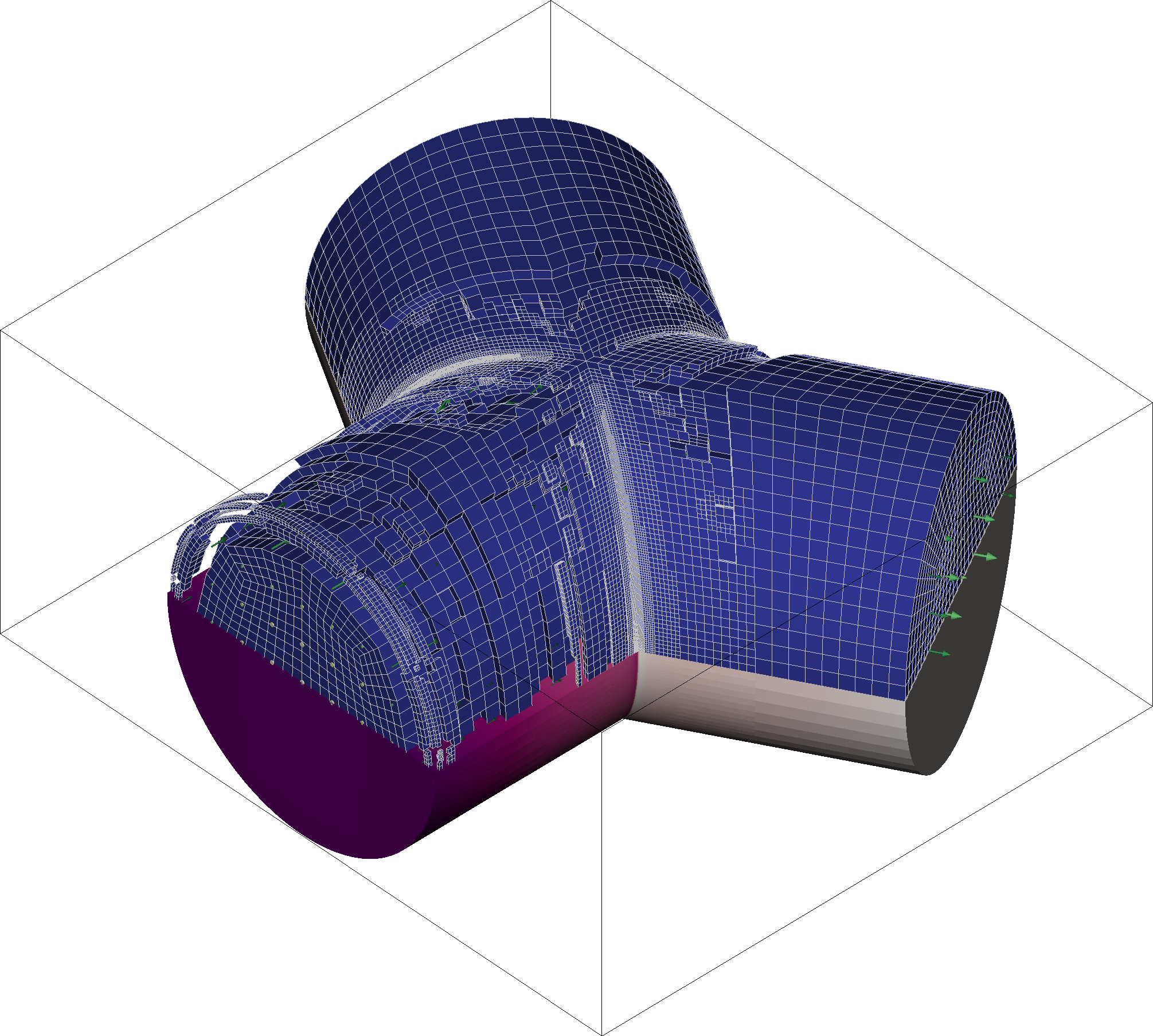}
  \end{minipage}
  \hfill
  \begin{minipage}{0.49\textwidth}
    \includegraphics[width=\textwidth]{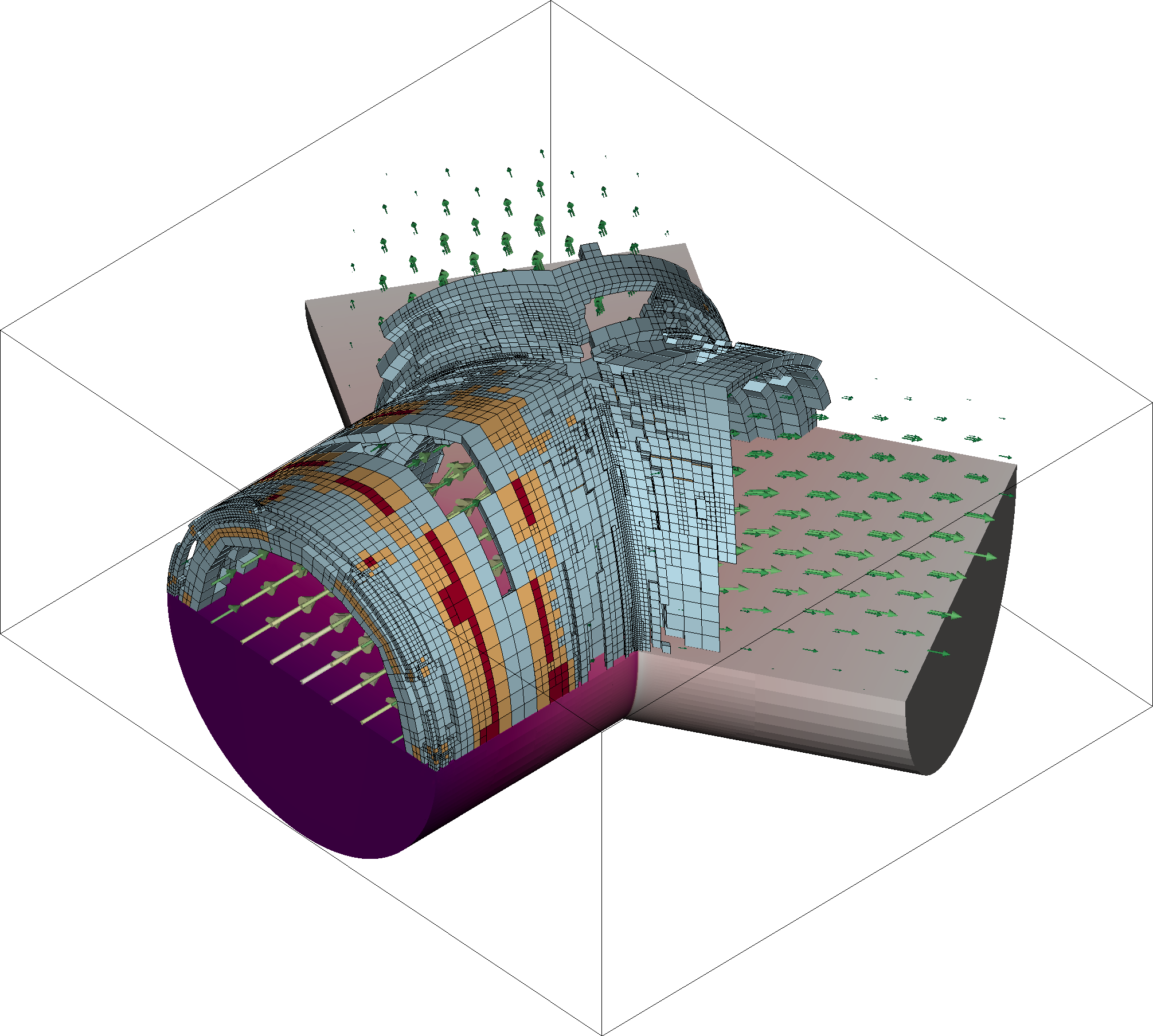}
  \end{minipage}
  \caption{\it Stokes flow through the Y-pipe as described by equation \eqref{eqs:stokes} after four adaptation cycles and three initial refinements.
  Top: The lower half shows the domain and the pressure, while the upper half depicts the mesh and a vector plot describing the velocity field.
  Bottom left: A cut-away showing those cells with low polynomial degrees, located generally where either the estimated errors are low or near the non-convex parts of the domain.
  Bottom right: Cut-away showing those cells with a high polynomial degree.
  %In the pictures in the bottom row, color represents the polynomial degree.
  }
  \label{fig:stokes_approximation}
\end{figure}

As a second test case, we consider the solution of the Stokes equation describing slow flow,
\begin{subequations}
\label{eqs:stokes}
\begin{align}
- \Delta \vec{u} + \nabla p &= \vec{0}, \\
- \nabla \cdot \vec{u} &= 0.
\end{align}%
\end{subequations}%
As domain, we choose a forked, ``Y''-shaped pipe, see Fig.~\ref{fig:stokes_approximation}. We impose no-slip boundary conditions on the lateral surfaces ($\vec{u} = 0$), and model the inflow at one opening as a Poisseuille flow via Dirichlet boundary conditions. The other two ends are modeled via zero-traction boundary conditions. Velocity and pressure solutions are also shown in Fig.~\ref{fig:stokes_approximation}.

The ``welding seams'' at which the three pipes meet are non-convex parts of the boundary, again resulting in singular solutions where we expect that the gradient of the velocity $\vec u$ becomes infinite; the pressure is also singular at these locations. We chose this as the second test case because it enables us to verify that enumerating \glspl{dof}, along with all of the other ingredients of our \hp-adaptive solution approach, are efficient and scale well also for three-dimensional problems with the much more complex choice of finite element and solver techniques necessary to solve the Stokes problem.

In particular, we use ``Taylor-Hood'' type elements $\pmb{Q}_k/Q_{k-1}$ \cite{taylor1973}, where the three components of the velocity solution use elements of polynomial degree $k$ and the single component of the pressure uses an element of polynomial degree ($k-1$). In our study, we choose a collection of elements with $k = 3,\dots,6$.

Both refinement and \hp-decision indicators are based on the scalar-valued pressure solution. We mark
\SI{10}{\percent} of cells for refinement and \SI{1}{\percent} for coarsening, which we divide equally into being \h- and \p-adapted.

We solve the linear saddle point system that results from discretization using flexible \gmres{} \cite{saad1993} and a Silvester-Wathen-type preconditioner \cite{silvester1994} in which we treat the elliptic block with the \ml{} \gls{amg} preconditioner \cite{gee2007} of \trilinos{} \cite{heroux2005}. This combination of solver and preconditioner is known to scale to very large problems, at least for elements of fixed order, see \cite{bangerth2012, kronbichler2012a}.

\begin{remark}
In our experiments, we have found that the \gls{amg} solver used in both of the test cases struggles with increasing fragmentation of polynomial degrees in the mesh. In order to address this, we limit the difference of polynomial degrees on neighboring cells to one, in a scheme not dissimilar to the commonly used approach of only allowing neighboring cells to differ by at most one level in mesh refinement. In our experiments, this ``smoothing'' of polynomial degrees reduces the number of solver iterations by up to \SI{70}{\percent}; this translates equally to the wallclock time spent on solving the linear system.
\end{remark}

\subsection{Load balancing}
\label{sec:results-load-balancing}

As mentioned in Section~\ref{sec:load_balancing}, it is not clear {\it a priori} how to weigh the contribution of each cell of a mesh to the overall cost of a program. As a consequence -- and unlike the \h-adaptive case -- it is not clear what the optimal load balancing strategy is.

Using the weighting proposed in Section~\ref{sec:load_balancing}, we have therefore run numerical experiments that vary the relative weighting of cells based on the number of \glspl{dof} on each cell. We carry out investigations on a mesh with a wide variety of polynomial degrees and a substantial number of hanging nodes that we obtain after a number of mesh refinement cycles. We keep this particular mesh, but partition it differently onto the available \gls{mpi} processes for varying values of the weighting exponent $c$ in \eqref{eq:load-balance-weight}, and run a complete refinement cycle involving enumeration of \glspl{dof}, assembly of the linear system, and solution of the linear system on the so-partitioned mesh.
These experiments were run on a workstation with two AMD EPYC\texttrademark{} 7552 processors (with 48 cores each, running at \SI{2.2}{\giga\hertz}) and \SI{512}{\giga\byte} of memory.

\begin{figure}
  \begin{tikzpicture}
\begin{groupplot}[
  scale=0.8,
  group style={
    group size=2 by 2,
    vertical sep=0pt
  },
%   width=9.5cm,
%   xlabel={Weighting exponent $c$},
%   ylabel={Wallclock time [\si{\second}] / Execution},
%   legend cell align=left,
%   legend pos=outer north east,
%   legend style={font=\tiny}
  ]

  % laplace top
  \nextgroupplot[
    axis y discontinuity = parallel,
    xticklabels = {},
    axis x line* = top,
    title = {Laplace},
    legend to name = legend:weights,
    legend columns=4,
    legend cell align=left,
    legend style={anchor=south, /tikz/every even column/.append style={column sep=0.5cm}},
    ymin=174,
    ymax=234,
    grid=none
  ]
  
  \addplot+ [only marks] table [y=solve, x=weighting_exponent, col sep=comma] {data/laplace/weights.csv};
  \addlegendentry{linear solver};
  
  % mimic 'assemble linear system' plot
  \addlegendimage{only marks, index of colormap=1 of Dark2, mark=square*, style={solid, fill=\pgfplotsmarklistfill}};
  \addlegendentry{assemble linear system};
  
  % \pgfplotsset{cycle list shift=1};
  \addplot+ [only marks, index of colormap=7 of Dark2, mark=star] table [y=full_cycle, x=weighting_exponent, col sep=comma] {data/laplace/weights.csv};
  \addlegendentry{full cycle};
  
%  \addplot [no markers, black, ultra thick] gnuplot [raw gnuplot] {
%    f(x) = a*x**2+b*x+c;
%    a=5;b=-25;c=200;
%    set datafile separator comma;
%    fit f(x) 'data/laplace/weights.csv' using 1:18 via a,b,c; % column 26 is solve, 18 is full_cycle
%    plot [x=1:4] f(x);
%    set print "parameters_laplace.dat";
%    print a, b, c;
%  };
  % gnuplot requires shell-escape, which is not allowed on arxiv.
  % instead, use fitted parameters directly.
   \addplot [no markers, black, ultra thick, domain=1:4] {8.56265731942037*\x^2 - 32.5839317129668*\x + 223.02245001861};
  \addlegendentry{quadratic fit};
  % \pgfplotsset{cycle list shift=-1};

  % stokes top
  \nextgroupplot[
    axis y discontinuity = parallel,
    xticklabels = {},
    axis x line*=top,
    title = {Stokes\vphantom{p}},
    ymin = 500,
    ymax = 1100,
    grid=none
  ]
  
  \addplot+ [only marks] table [y=solve, x=weighting_exponent, col sep=comma] {data/stokes/weights.csv};
  
  % \pgfplotsset{cycle list shift=1};
  \addplot+ [only marks, index of colormap=7 of Dark2, mark=star] table [y=full_cycle, x=weighting_exponent, col sep=comma] {data/stokes/weights.csv};
  
%  \addplot [no markers, black, ultra thick] gnuplot [raw gnuplot] {
%    f(x) = a*x**2+b*x+c;
%    a=40;b=-200;c=1100;
%    set datafile separator comma;
%    fit f(x) 'data/stokes/weights.csv' using 1:15 via a,b,c; % column 23 is solve, 15 is full_cycle
%    plot [x=1:4] f(x);
%    set print "parameters_stokes.dat";
%    print a, b, c;
%  };
  % gnuplot requires shell-escape, which is not allowed on arxiv.
  % instead, use fitted parameters directly.
   \addplot [no markers, black, ultra thick, domain=1:4] {43.4851032396886*\x^2 - 211.090556510164*\x + 1121.88802633854};

  % laplace bottom
  \nextgroupplot[
    axis x line* = bottom,
    y post scale = 0.27,
    ymin = 0,
    ymax = 20,
    ytick distance = 20,
    grid=none
  ]
  
  \pgfplotsset{cycle list shift=1};
  
  \addplot+ [only marks] table [y=assemble_system, x=weighting_exponent, col sep=comma] {data/laplace/weights.csv};

  % stokes bottom
  \nextgroupplot[
    axis x line*=bottom,
    y post scale = 0.27,
    ymin = 0,
    ymax = 200,
    ytick distance = 200,
    grid=none
  ]
  
  \pgfplotsset{cycle list shift=1};
  
  \addplot+ [only marks] table [y=assemble_system, x=weighting_exponent, col sep=comma] {data/stokes/weights.csv};
  \end{groupplot}

  % shared xlabel
  \path (group c1r2.south east) -- node[below=3.5ex]{Weighting exponent $c$} (group c2r2.south west);
  
  % shared ylabel
  \path (group c1r1.north west) -- node[left=3em, anchor=center, rotate=90]{Wallclock time [\si{s}] / Execution} (group c1r2.south west);
  
  % shared legend
  \path (group c1r1.north east) -- node[above=5ex]{\ref*{legend:weights}} (group c2r1.north west);
\end{tikzpicture}%
  \caption[Wall times for load balancing with varying weighting exponents.]{\it Wall clock times for several operations of a complete adaptation cycle, when partitioning the mesh using different weighting exponents $c$, see \eqref{eq:load-balance-weight}.
  Left: For one cycle of the two-dimensional Laplace problem of Section~\ref{sec:test-case-1}. The problem has about 51 million \glspl{dof} and is solved on 96 \gls{mpi} processes.
  Right: For one cycle of the three-dimensional Stokes problem of Section~\ref{sec:test-case-2}. The problem has about 15 million \glspl{dof} and is solved on 96 \gls{mpi} processes.}
  \label{fig:results_weights}
\end{figure}
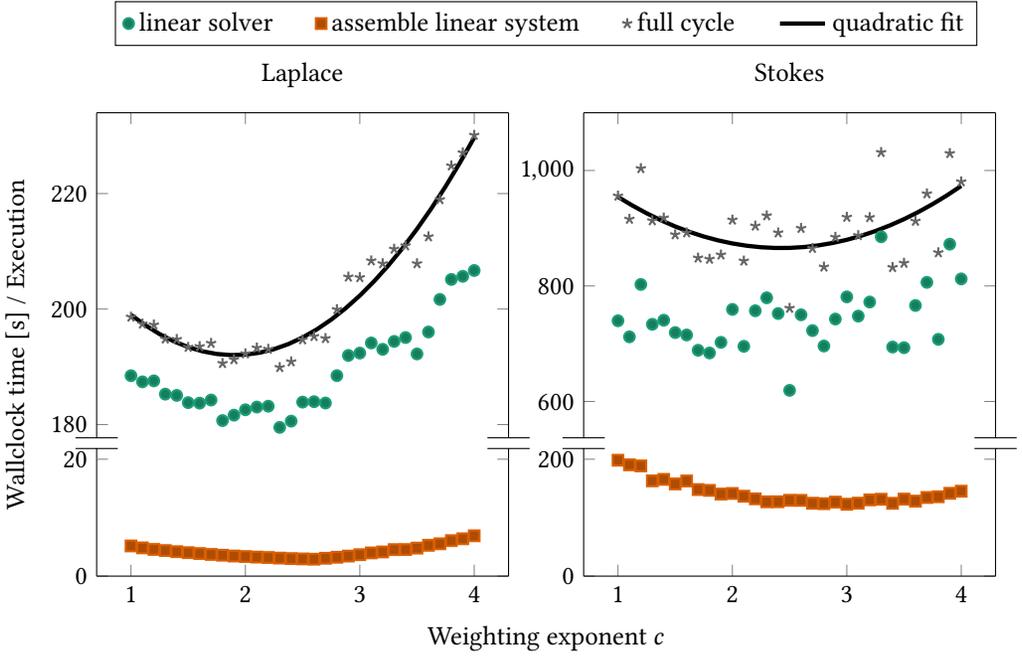

The results are shown in Fig.~\ref{fig:results_weights} for the two
test cases defined in Section~\ref{sec:testcases}. For both cases, the
largest contribution to the overall cost is the linear solver; the
noise in the corresponding curves results from slightly different
numbers of linear solver iterations, likely a consequence of decisions
made in how the \gls{amg} algorithm builds its hierarchy in response
to which rows of the overall matrix are stored on which process. The
data shown in Fig.~\ref{fig:results_weights} suggest that the overall
run time is minimized with an exponent of $c\approx 1.9$ for the
Laplace test case, and $c\approx 2.4$ for the Stokes test case. We use these values for the weighting exponent for all other
experiments shown below.

The data shown in the figure makes it clear
that the optimal exponent depends on the problem solved, and needs to
be assessed for each problem individually. However, in general the
dependency of the run time on the specific choice of exponent is relatively weak.

\subsection{Efficiency and scalability of algorithms}
\label{sec:results-scaling}

In the following, we assess whether the algorithms we proposed in Section~\ref{sec:enumeration} are efficient and scale to large problem sizes. To answer this question, we first discuss what ``scalability'' means in the context of \hp-adaptive methods, before turning to results obtained on the test cases defined in Section~\ref{sec:testcases}.

All results shown in this subsection were obtained on the \expanse{} supercomputer. Each standard computing node is equipped with two AMD EPYC\texttrademark{} 7742 processors (with 64 cores each, running at \SI{2.25}{\giga\hertz}) and \SI{256}{\giga\byte} of DDR4 DRAM memory. Communication between nodes happens via a Mellanox\textsuperscript{\textregistered} HDR-100 InfiniBand Interconnect network operating in a hybrid fat-tree topology. More information on the configuration of the machine can be found in \cite{strande2021}.

\subsubsection{How to define scalability?}
\label{sec:results-scalability-definition}

One typically measures the efficiency of a parallel algorithm running on $P$ processes operating in parallel on $N$ work items through either ``strong scaling'' (where the problem size $N$ is fixed and we vary the number of processes $P$) or ``weak scaling'' (where one increases the problem size $N$ along with the number of processes $P$, keeping $N/P$ constant). In both cases, one measures the time it takes the algorithm to complete work.

For \h-adaptive algorithms, it is relatively straightforward to define what $N$ is supposed to be: it could be (i) the number of cells in the mesh, (ii) the number of unknowns in a finite element discretization on that mesh (which equals the size of the linear systems that result), or (iii) the number of nonzero entries in the matrix (which determines the cost of a matrix-vector product, but is also an important consideration in the cost of algorithms such as \gls{amg}). The choice of which of these we want to call $N$ is unimportant \textit{because they are all proportional to each other}. Indeed, if one uses an optimal solver such as multigrid, one could also (iv) define $N$ to be the number of floating point operations required to solve the linear system for a given problem -- it is again proportional to the other measures.

But things are not this easy for \hp-adaptive methods: when using different polynomial degrees on cells, the four quantities mentioned above are no longer proportional to each other when considering an \hp-fragmented mesh. This disproportionality is of no importance when considering strong scalability, because the problem size $N$ is fixed. But it is not obvious how to define weak scalability because a sequence of problems that keeps $N/P$ constant for one definition of $N$ may not imply that $N/P$ is constant for any of the other definitions of $N$. Similarly, we show results below where we increase $N$ for fixed $P$, observing how time scales with $N$ -- for which, again, the observed scaling depends on what definition of $N$ we choose.

As a consequence, we describe results below where we either use $N=N_\text{\glsfmtshortpl{dof}}$ (the number of global \glspl{dof} in the problem), or $N=N_\text{nonzeros}$ (the global number of nonzero entries in the matrix which needs to be solved with on a given mesh). As expected, we will see that operations such as the assembly of a linear system and its solution do not scale with the problem size as ${\mathcal O}(N_\text{\glsfmtshortpl{dof}})$, but they instead scale close to the amount of work as ${\mathcal O}(N_\text{nonzeros})$.

\subsubsection{Results for the Laplace test case of Section~\ref{sec:test-case-1}.}

\begin{figure}
  \begin{tikzpicture}
  \begin{groupplot}[
    group style = {
      group size = 2 by 1,
    %   horizontal sep = 0.1\textwidth,
    },
  ]
  
  \nextgroupplot[
    width = 0.52\textwidth,
    xmode = log,
    ymode = log,
    title = 1024 \gls{mpi} processes,
    ylabel = {Wallclock time [\si{\second}] / Execution},
    legend to name = legend:laplace-cycles,
    legend columns=3,
    legend cell align=left,
    legend style={anchor=south, /tikz/every even column/.append style={column sep=0.5cm}},
    % xmin = 1e7,
    % xmax = 2e9,
    ymin = 8e-3,
    ymax = 5e2,
    grid = none,
  ]
  
  % data
  \addplot table [y=solve, x=dofs, col sep=comma] {data/laplace/cycles_proc1024.csv};
  \addlegendentry{linear solver};

  \addplot table [y=assemble_system, x=dofs, col sep=comma] {data/laplace/cycles_proc1024.csv};
  \addlegendentry{assemble linear system};

  \addplot table [y=setup_system, x=dofs, col sep=comma] {data/laplace/cycles_proc1024.csv};
  \addlegendentry{setup data structures};
  
  \addplot table [y=distribute_dofs, x=dofs, col sep=comma]
  {data/laplace/cycles_proc1024.csv};
  \addlegendentry{enumerate \glspl{dof}};

  \addplot table [y=estimate_mark, x=dofs, col sep=comma] {data/laplace/cycles_proc1024.csv};
  \addlegendentry{estimate and mark};

  \addplot table [y=refine, x=dofs, col sep=comma] {data/laplace/cycles_proc1024.csv};
  \addlegendentry{coarsen and refine};
  
  % auxiliary lines
  \begin{scope}
    \draw[dashed, very thick] ({axis cs:102400000,0}|-{rel axis cs:0,1}) -- ({axis cs:102400000,0}|-{rel axis cs:0,0});
  \end{scope}
  \addlegendimage{dashed, very thick};
  \addlegendentry{$10^5$ \glspl{dof} per process};
  
%   \begin{scope}
%     \draw[dotted, very thick] ({axis cs:51887592,0}|-{rel axis cs:0,1}) -- ({axis cs:51887592,0}|-{rel axis cs:0,0});
%   \end{scope}
%   \addlegendimage{dotted, very thick};
%   \addlegendentry{all $k$ appear in domain};
  
  % optimal domain=1e7:2e9
  \addplot[thick, samples=2, domain=12591105:697832457] {2.5e-7*x};
  \addplot[thick, samples=2, domain=12591105:697832457] {0.5e-8*x};
  \addplot[thick, samples=2, domain=12591105:697832457] {0.5e-9*x};
  \addlegendentry{$\mathcal{O}(N_\text{\glsfmtshortpl{dof}})$};
  
  \nextgroupplot[
    width = 0.52\textwidth,
    xmode = log,
    ymode = log,
    % xlabel = {Number of \glspl{dof}},
    % xmin = 1e7,
    % xmax = 2e9,
    ymin = 8e-3,
    ymax = 5e2,
    grid = none,
    title = 4096 \gls{mpi} processes,
  ]
  
  % data
  \addplot table [y=solve, x=dofs, col sep=comma] {data/laplace/cycles_proc4096.csv};

  \addplot table [y=assemble_system, x=dofs, col sep=comma] {data/laplace/cycles_proc4096.csv};

  \addplot table [y=setup_system, x=dofs, col sep=comma] {data/laplace/cycles_proc4096.csv};
  
  \addplot table [y=distribute_dofs, x=dofs, col sep=comma]
  {data/laplace/cycles_proc4096.csv};

  \addplot table [y=estimate_mark, x=dofs, col sep=comma] {data/laplace/cycles_proc4096.csv};
 
  \addplot table [y=refine, x=dofs, col sep=comma] {data/laplace/cycles_proc4096.csv};

  % auxiliary lines
  \begin{scope}
    \draw[dashed, very thick] ({axis cs:409600000,0}|-{rel axis cs:0,1}) -- ({axis cs:409600000,0}|-{rel axis cs:0,0});
  \end{scope}
  
%   \begin{scope}
%     \draw[dotted, very thick] ({axis cs:53537286,0}|-{rel axis cs:0,1}) -- ({axis cs:53537286,0}|-{rel axis cs:0,0});
%   \end{scope}
  
  % optimal domain=1e7:2e9
  \addplot[thick, samples=2, domain=12591105:1445633269] {2.5e-7*x/4};
  \addplot[thick, samples=2, domain=12591105:1445633269] {0.5e-8*x/4};
  \addplot[thick, samples=2, domain=12591105:1445633269] {0.5e-9*x/4};
\end{groupplot}
  
  % shared xlabel
  \path (group c1r1.south east) -- node[below=3.5ex]{Number of \glspl{dof}} (group c2r1.south west);
  
  % shared legend
  \path (group c1r1.north east) -- node[above=5ex]{\ref*{legend:laplace-cycles}} (group c2r1.north west);
\end{tikzpicture}%
  \caption{\it Laplace problem: Scaling of wallclock time as a function of the number of unknowns $N_\text{\glsfmtshortpl{dof}}$ on a sequence of consecutively refined meshes, for \num{1024} (left) and \num{4096} \gls{mpi} processes. Each \gls{mpi} process owns more than $10^5$ \glspl{dof} only to the right of the indicated vertical line; to the left of this line, processes do not have enough work to offset the cost of communication, and parallel efficiency should not be expected.
%   Each finite element is represented at least once in the mesh only to the right side of the designated vertical line.
  The solid black trend lines indicate optimal scaling, ${\mathcal
    O}(N_\text{\glsfmtshortpl{dof}})$.
\\
Since the computations on the right were done on four times as many
processes as on the left, we have offset the black trend lines
downward by a factor of four: Assuming optimal strong scaling when
increasing the number of \gls{mpi} processes, a point located on the
trend line in the left sub-figure should have a corresponding point on
the offset trend line also in the right sub-figure.
  }
  \label{fig:laplace_cycles}
\end{figure}
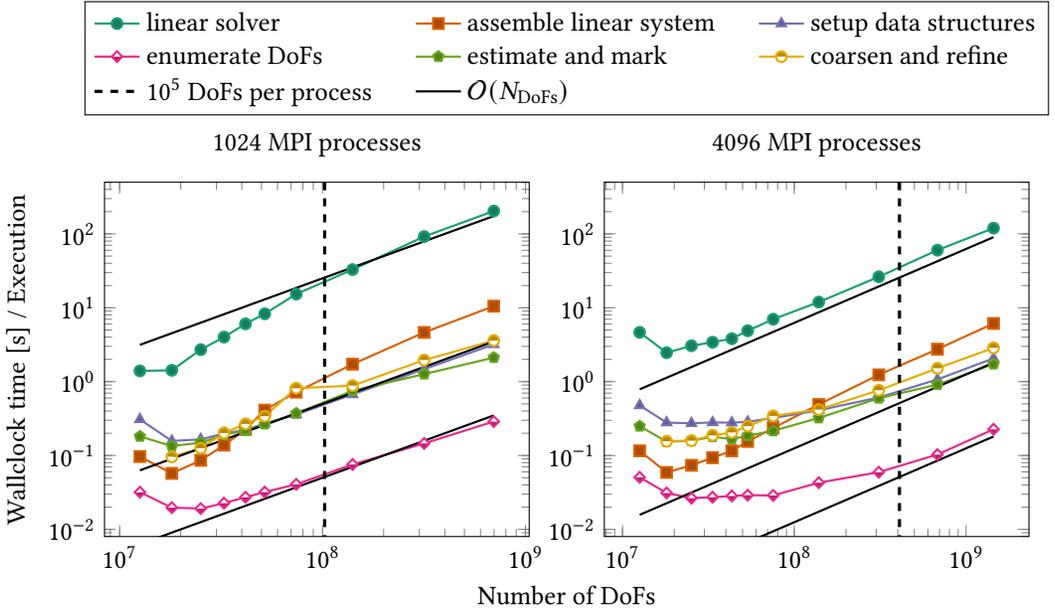

With these considerations in mind, let us now turn to concrete timing data. Below, we show results for how much time our implementation of the Laplace test case of Section~\ref{sec:test-case-1} spends in each of the following categories of operations (ordered roughly in their relevance to the overall run time to the program):

\begin{itemize}
  \item {\it Linear solver:} This category includes setting up the \gls{amg} preconditioner, and then solving the linear system.
  
  \item {\it Assemble linear system:}
  Compute cell-local matrix and right-hand side vector contributions
  to the linear system, and insertion into the global objects. This
  step also includes communicating these contributions to the process owning a matrix or vector row if necessary.
  
  \item {\it Setup data structures:}
  This step includes a number of setup steps that happen after generating a mesh and before the assembly of the linear system. Specifically, we include the enumeration of \glspl{dof};
  exchanging between processes which non-locally owned matrix entries they will write into;
  setting up a sparsity pattern for the global matrix;
  allocation of memory for the system matrix and vectors; and
  determining constraints that result from hanging nodes and boundary conditions.
  
  \item {\it Enumerate \glspl{dof}:}
  This category, a subset of the previous one, measures the time to
  enumerate all \glspl{dof} based on the algorithm discussed in Section~\ref{sec:enum-hp-parallel}.

  \item {\it Estimate and mark:}
  Once the linear system has been solved, this step computes error and smoothness estimates for each locally owned cell. It then
  computes global thresholds for \hp-adaptation, and flags cells for either \h- or \p-adaptation.
  
  \item {\it Coarsen and refine:}
  This final step performs the actual \h-adaptation on marked cells while enforcing a 2:1 cell size relationship across faces. It also
  updates the associated finite element on cells (\p-adaptation) while limiting the difference of polynomial degrees across cell interfaces.
  This category also measures the transfer data between old and new mesh, as well as the cost of re-partitioning the mesh between processes.
\end{itemize}

Fig.~\ref{fig:laplace_cycles} shows timing information for a situation where we repeatedly solve the problem while adaptively refining the \hp-mesh, on both \num{1024} and \num{4096} MPI processes. In this setup, with a fixed number $P$ of processes, one would hope that the run time increases linearly with the problem size $N$.
% The figure shows that that is largely true where we have here chosen $N=N_\text{\glspl{dof}}$ as the number of \glspl{dof} on each of the meshes 
Our results demonstrate that this linearity holds when $N$ is the $N_\text{\glspl{dof}}$ on each of the meshes -- at least once the problem is large enough.
Importantly for the current paper, operations such as estimating \hp-indicators and refining the mesh accordingly, and in particular the enumeration of \glspl{dof} using the algorithm of Section~\ref{sec:enum-hp-parallel} are only minor contributions to the overall run time, which is dominated by the assembly and in particular solution of linear systems.

\begin{figure}
  \begin{tikzpicture}
  \begin{groupplot}[
    group style = {
      group size = 2 by 1,
    %   horizontal sep = 0.1\textwidth,
    },
  ]
  
  \nextgroupplot[
    width = 0.52\textwidth,
    xmode = log,
    ymode = log,
    % xmin = 1e8,
    % xmax = 1e11,
    ymin = 8e-3,
    ymax = 5e2,
    title = 1024 \gls{mpi} processes,
    ylabel = {Wallclock time [\si{\second}] / Execution},
    legend to name = legend:laplace-cycles-nonzeros,
    legend columns=3,
    legend cell align=left,
    grid = none,
    legend style={anchor=south, /tikz/every even column/.append style={column sep=0.5cm}},
  ]
  
  % data
  \addplot table [y=solve, x=nonzero_elements, col sep=comma] {data/laplace/cycles_proc1024.csv};
  \addlegendentry{linear solver};

  \addplot table [y=assemble_system, x=nonzero_elements, col sep=comma] {data/laplace/cycles_proc1024.csv};
  \addlegendentry{assemble linear system};

  \addplot table [y=setup_system, x=nonzero_elements, col sep=comma] {data/laplace/cycles_proc1024.csv};
  \addlegendentry{setup data structures};
  
  \addplot table [y=distribute_dofs, x=nonzero_elements, col sep=comma]
  {data/laplace/cycles_proc1024.csv};
  \addlegendentry{enumerate \glspl{dof}};

  \addplot table [y=estimate_mark, x=nonzero_elements, col sep=comma] {data/laplace/cycles_proc1024.csv};
  \addlegendentry{estimate and mark};

  \addplot table [y=refine, x=nonzero_elements, col sep=comma] {data/laplace/cycles_proc1024.csv};
  \addlegendentry{coarsen and refine};
  
  % auxiliary lines
  \begin{scope}
    \draw[dashed, very thick] ({axis cs:1024000000,0}|-{rel axis cs:0,1}) -- ({axis cs:1024000000,0}|-{rel axis cs:0,0});
  \end{scope}
  \addlegendimage{dashed, very thick};
  \addlegendentry{$10^6$ nonzeros per process};
  
%   \begin{scope}
%     \draw[dotted, very thick] ({axis cs:1694948604,0}|-{rel axis cs:0,1}) -- ({axis cs:1694948604,0}|-{rel axis cs:0,0});
%   \end{scope}
%   \addlegendimage{dotted, very thick};
%   \addlegendentry{all $k$ appear in domain};
  
  % optimal
  \addplot[thick, samples=2, domain=201048153:50861845017] {10^(-8.4)*x};
  \addplot[thick, samples=2, domain=201048153:50861845017] {10^(-9.7)*x};
  \addlegendentry{$\mathcal{O}(N_\text{nonzeros})$};
  
  \nextgroupplot[
    width = 0.52\textwidth,
    xmode = log,
    ymode = log,
    grid = none,
    % xmin = 1e8,
    % xmax = 2e11,
    ymin = 8e-3,
    ymax = 5e2,
    title = 4096 \gls{mpi} processes,
  ]
  
  % data
  \addplot table [y=solve, x=nonzero_elements, col sep=comma] {data/laplace/cycles_proc4096.csv};

  \addplot table [y=assemble_system, x=nonzero_elements, col sep=comma] {data/laplace/cycles_proc4096.csv};

  \addplot table [y=setup_system, x=nonzero_elements, col sep=comma] {data/laplace/cycles_proc4096.csv};
  
  \addplot table [y=distribute_dofs, x=nonzero_elements, col sep=comma]
  {data/laplace/cycles_proc4096.csv};

  \addplot table [y=estimate_mark, x=nonzero_elements, col sep=comma] {data/laplace/cycles_proc4096.csv};
 
  \addplot table [y=refine, x=nonzero_elements, col sep=comma] {data/laplace/cycles_proc4096.csv};

  % auxiliary lines
  \begin{scope}
    \draw[dashed, very thick] ({axis cs:4096000000
,0}|-{rel axis cs:0,1}) -- ({axis cs:4096000000
,0}|-{rel axis cs:0,0});
  \end{scope}
  
%   \begin{scope}
%     \draw[dotted, very thick] ({axis cs:1757526006,0}|-{rel axis cs:0,1}) -- ({axis cs:1757526006,0}|-{rel axis cs:0,0});
%   \end{scope}
  
  % optimal
  \addplot[thick, samples=2, domain=201048153:108532301149] {10^(-8.4)*x/4};
  \addplot[thick, samples=2, domain=201048153:108532301149] {10^(-9.7)*x/4};
  \end{groupplot}
  
  % shared xlabel
  \path (group c1r1.south east) -- node[below=3.5ex]{Number of nonzero matrix elements} (group c2r1.south west);
  
  % shared legend
  \path (group c1r1.north east) -- node[above=5ex]{\ref*{legend:laplace-cycles-nonzeros}} (group c2r1.north west);
\end{tikzpicture}%
  \caption{\it Laplace problem: The same scaling data as shown in Fig.~\ref{fig:laplace_cycles}, except shown as a function of the nonzero entries of the matrix.}
  \label{fig:laplace_cycles_nonzeros}
\end{figure}
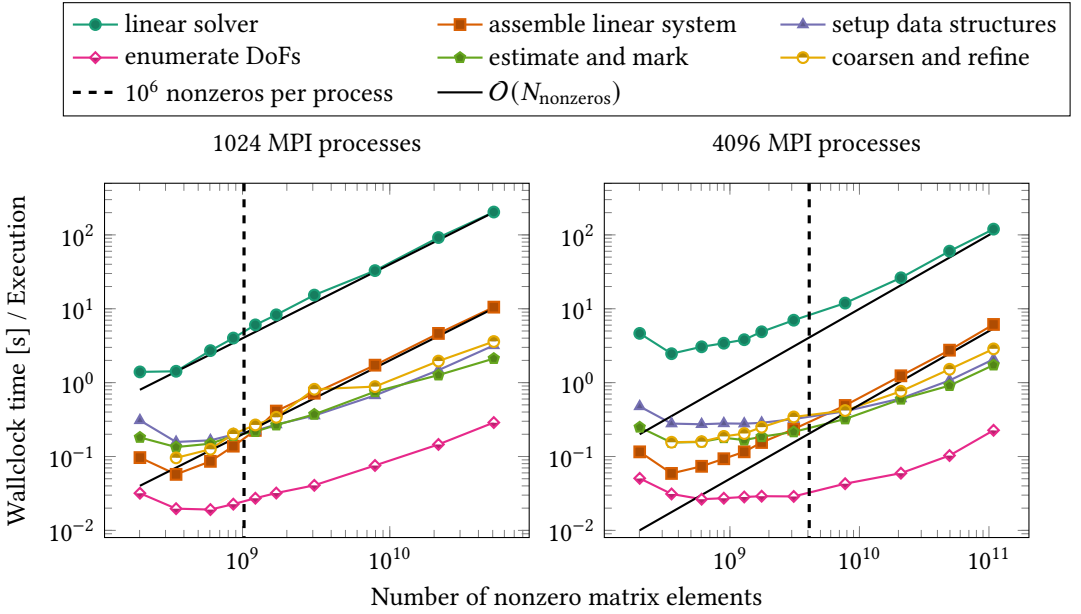

On the other hand, Fig.~\ref{fig:laplace_cycles} also shows that both the assembly and the solution of the linear system do not scale like $\mathcal{O}(N_\text{\glspl{dof}})$. This result may not be surprising in view of the discussions of Section~\ref{sec:results-scalability-definition}: as we move from left to right, we do not only increase the number of unknowns, but also increase polynomial degrees on cells, resulting in denser and denser linear systems that are more costly to assemble and solve. As a consequence, Fig.~\ref{fig:laplace_cycles_nonzeros} shows the same data as a function of the nonzero entries $N_\text{nonzeros}$. This figure illustrates that using this definition, both assembly and the solution of linear systems scale nearly perfectly as $\mathcal{O}(N_\text{nonzeros})$.

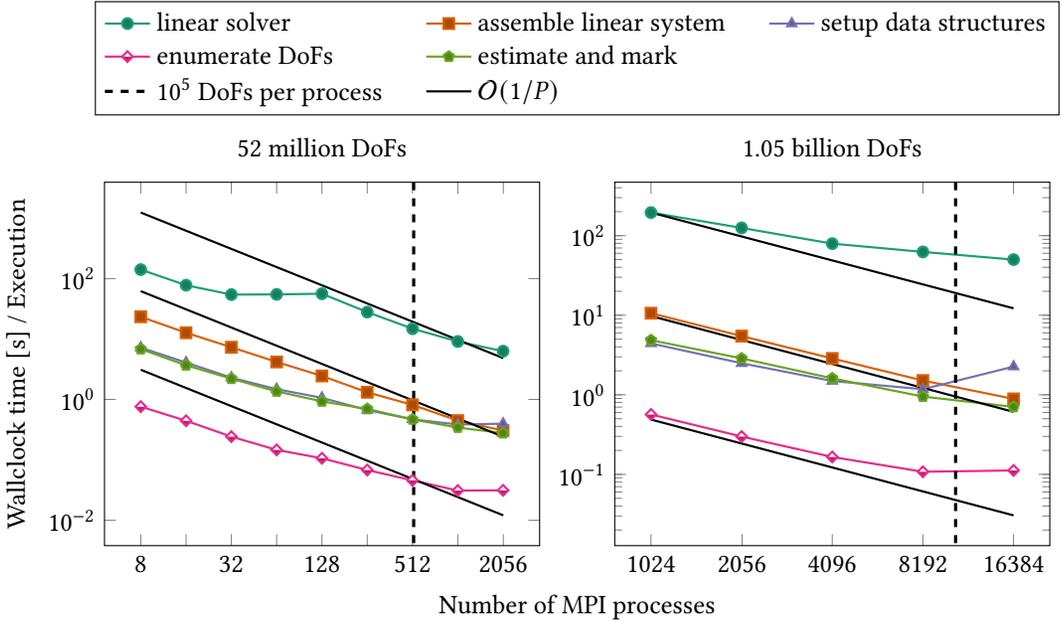
\begin{figure}
  \begin{tikzpicture}
  \begin{groupplot}[
    group style = {
      group size = 2 by 1,
    %   horizontal sep = 0.1\textwidth,
    },
  ]
  
  \nextgroupplot[
    width = 0.53\textwidth,
    xmode = log,
    ymode = log,
    title = 52 million \glspl{dof},
    ylabel = {Wallclock time [\si{\second}] / Execution},
    xtick = {8, 16, 32, 64, 128, 256, 512, 1024, 2056},
    xticklabels = {8, , 32, , 128, , 512, , 2056},
    legend to name = legend:laplace-strong,
    legend columns=3,
    legend cell align=left,
    legend style={anchor=south, /tikz/every even column/.append style={column sep=0.5cm}},
    grid = none
  ]
  
  % data
  \addplot table [y=solve, x=processes, col sep=comma] {data/laplace/strong_nref10.csv};
  \addlegendentry{linear solver};

  \addplot table [y=assemble_system, x=processes, col sep=comma] {data/laplace/strong_nref10.csv};
  \addlegendentry{assemble linear system};

  \addplot table [y=setup_system, x=processes, col sep=comma] {data/laplace/strong_nref10.csv};
  \addlegendentry{setup data structures};
  
  \addplot table [y=distribute_dofs, x=processes, col sep=comma]
  {data/laplace/strong_nref10.csv};
  \addlegendentry{enumerate \glspl{dof}};

  \addplot table [y=estimate_mark, x=processes, col sep=comma] {data/laplace/strong_nref10.csv};
  \addlegendentry{estimate and mark};
  
  % empty legend entry
  \addlegendimage{empty legend};
  \addlegendentry{};
  
  % auxiliary lines
  \begin{scope}
    \draw[dashed, very thick] ({axis cs:521.25652,0}|-{rel axis cs:0,1}) -- ({axis cs:521.25652,0}|-{rel axis cs:0,0});
  \end{scope}
  \addlegendimage{dashed, very thick};
  \addlegendentry{$10^5$ \glspl{dof} per process};

  % optimal
  \addplot[thick, samples=2, domain=8:2056] {4.9551e-2*2e5/x};
  \addplot[thick, samples=2, domain=8:2056] {4.9551e-2*1e4/x};
  \addplot[thick, samples=2, domain=8:2056] {4.9551e-2*.5e3/x};
  \addlegendentry{$\mathcal{O}(1/P)$};
  
  \nextgroupplot[
    width = 0.53\textwidth,
    xmode = log,
    ymode = log,
    grid = none,
    title = 1.05 billion \glspl{dof},
    xtick = {1024, 2056, 4096, 8192, 16384},
    xticklabels = {1024, 2056, 4096, 8192, 16384},
  ]
  
  % data
  \addplot table [y=solve, x=processes, col sep=comma] {data/laplace/strong_nref12.csv};

  \addplot table [y=assemble_system, x=processes, col sep=comma] {data/laplace/strong_nref12.csv};

  \addplot table [y=setup_system, x=processes, col sep=comma] {data/laplace/strong_nref12.csv};
  
  \addplot table [y=distribute_dofs, x=processes, col sep=comma]
  {data/laplace/strong_nref12.csv};

  \addplot table [y=estimate_mark, x=processes, col sep=comma] {data/laplace/strong_nref12.csv};
  
  % auxiliary lines
  \begin{scope}
    \draw[dashed, very thick] ({axis cs:10519.63562,0}|-{rel axis cs:0,1}) -- ({axis cs:10519.63562,0}|-{rel axis cs:0,0});
  \end{scope}
  
  % optimal
  \addplot[thick, samples=2, domain=1024:16384] {2e5/x};
  \addplot[thick, samples=2, domain=1024:16384] {1e4/x};
  \addplot[thick, samples=2, domain=1024:16384] {.5e3/x};
  \end{groupplot}
  
  % shared xlabel
  \path (group c1r1.south east) -- node[below=3.5ex]{Number of \gls{mpi} processes} (group c2r1.south west);
  
  % shared legend from 
  \path (group c1r1.north east) -- node[above=5ex]{\ref*{legend:laplace-strong}} (group c2r1.north west);
\end{tikzpicture}%
  \caption{\it Laplace problem: Strong scaling for one advanced adaptation cycle at different problem sizes. Each \gls{mpi} process owns more than $10^5$ \glspl{dof} only to the left of the indicated vertical line; to the right of this line, processes do not have enough work to offset the cost of communication, and parallel efficiency should not be expected. Left: Fixed problem size of roughly 52 million \glspl{dof}. Right: Fixed problem size of roughly 1.05 billion \glspl{dof}. The trend lines for ${\mathcal O}(1/P)$ are offset between the two panels by the ratio of the size of the problem to allow for assessing weak scalability of the algorithms.}
  \label{fig:laplace_strong}
\end{figure}

A comparison of the left and right panels of Figs.~\ref{fig:laplace_cycles} and \ref{fig:laplace_cycles_nonzeros} -- and in particular how the various curves approach the trend lines $\mathcal{O}(N)$ that are offset in the panels by the ratio of the number of MPI processes used -- shows that for sufficiently large problems, we also have good \textit{strong} scaling. We expand on this in Fig.~\ref{fig:laplace_strong} where we show scaling for a fixed problem with the number of processes. The figure shows that most operations may not scale perfectly as $\mathcal{O}(1/P)$, but scalability is at least adequate as long as the problem size per process remains sufficiently large (to the left of the dashed line). The exception is the performance of the linear solver; this is a known problem with implementations of algebraic multigrid methods, but also beyond the scope of the current paper.

\subsubsection{Results for the Stokes test case of Section~\ref{sec:test-case-2}.}
We repeat many of these timing studies, using the same timing categories, for the Stokes test case to assess whether our results also hold for a more complex, three-dimensional, and vector-valued problem.

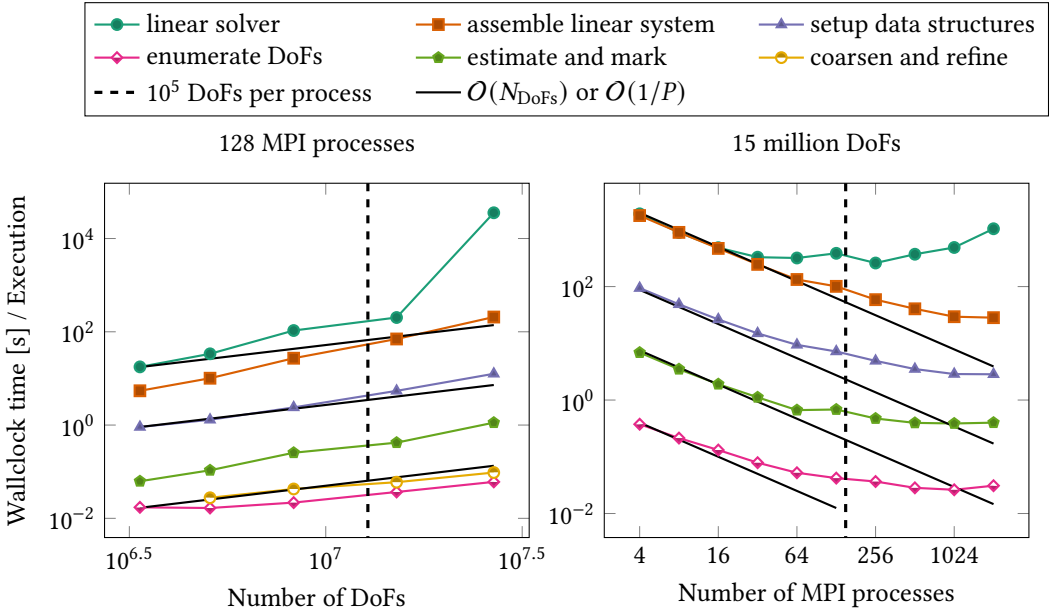
\begin{figure}
  \begin{tikzpicture}
  \begin{groupplot}[
    group style = {
      group size = 2 by 2,
    %   horizontal sep = 0.1\textwidth,
    },
  ]
  
  \nextgroupplot[
    width = 0.52\textwidth,
    xmode = log,
    ymode = log,
    % xmin = 2e6,
    % xmax = 2e7,
    % domain = 2e6:2e7,
    % ymin = 1e-2,
    % ymax = 3e3,
    title = 128 \gls{mpi} processes,
    xlabel = {Number of \glspl{dof}},
    ylabel = {Wallclock time [\si{\second}] / Execution},
    legend to name = legend:stokes-scaling,
    legend columns=3,
    legend cell align=left,
    legend style={anchor=south, /tikz/every even column/.append style={column sep=0.5cm}},
    grid = none
  ]
  
  % data
  \addplot table [y=solve, x=dofs, col sep=comma] {data/stokes/cycles_proc128.csv};
  \addlegendentry{linear solver};

  \addplot table [y=assemble_system, x=dofs, col sep=comma] {data/stokes/cycles_proc128.csv};
  \addlegendentry{assemble linear system};

  \addplot table [y=setup_system, x=dofs, col sep=comma] {data/stokes/cycles_proc128.csv};
  \addlegendentry{setup data structures};
  
  \addplot table [y=distribute_dofs, x=dofs, col sep=comma]
  {data/stokes/cycles_proc128.csv};
  \addlegendentry{enumerate \glspl{dof}};

  \addplot table [y=estimate_mark, x=dofs, col sep=comma] {data/stokes/cycles_proc128.csv};
  \addlegendentry{estimate and mark};

  \addplot table [y=refine, x=dofs, col sep=comma] {data/stokes/cycles_proc128.csv};
  \addlegendentry{coarsen and refine};

  % auxiliary lines
  \begin{scope}
    \draw[dashed, very thick] ({axis cs:12800000,0}|-{rel axis cs:0,1}) -- ({axis cs:12800000,0}|-{rel axis cs:0,0});
  \end{scope}
  \addlegendimage{dashed, very thick};
  \addlegendentry{$10^5$ \glspl{dof} per process};
  
  % optimal
  % \addplot[thick, samples=2, domain=3365108:26777965] {1e-5*x};
  % \addplot[thick, samples=2, domain=3365108:26777965] {4.5e-7*x};
  \addplot[thick, samples=2, domain=3365108:26777965] {5.2e-6*x};
  \addplot[thick, samples=2, domain=3365108:26777965] {2.7e-7*x};
  \addplot[thick, samples=2, domain=3365108:26777965] {5e-9*x};
  \addlegendentry{$\mathcal{O}(N_\text{\glsfmtshortpl{dof}})$ or $\mathcal{O}(1/P)$};
  
  \nextgroupplot[
    width = 0.52\textwidth,
    xmode = log,
    ymode = log,
    % xmin = 2,
    % xmax = 4096,
    % ymin = 1e-2,
    % ymax = 3e3,
    title = {15 million \glspl{dof}\vphantom{p}},
    xlabel = {Number of \gls{mpi} processes},
    xtick = {4, 16, 64, 256, 1024},
    xticklabels = {4, 16, 64, 256, 1024},
    grid = none
  ]
  
    % data
  \addplot table [y=solve, x=processes, col sep=comma] {data/stokes/strong_nref3.csv};

  \addplot table [y=assemble_system, x=processes, col sep=comma] {data/stokes/strong_nref3.csv};

  \addplot table [y=setup_system, x=processes, col sep=comma] {data/stokes/strong_nref3.csv};
  
  \addplot table [y=distribute_dofs, x=processes, col sep=comma]
  {data/stokes/strong_nref3.csv};

  \addplot table [y=estimate_mark, x=processes, col sep=comma] {data/stokes/strong_nref3.csv};

  % auxiliary lines
  \begin{scope}
    \draw[dashed, very thick] ({axis cs:151.10817,0}|-{rel axis cs:0,1}) -- ({axis cs:151.10817,0}|-{rel axis cs:0,0});
  \end{scope}
  
  % optimal
  \addplot[thick, samples=2, domain=4:2048] {1e3*8/x};
  \addplot[thick, samples=2, domain=4:2048] {3.5e2/x};
  \addplot[thick, samples=2, domain=4:2048] {3e1/x};
  \addplot[thick, samples=2, domain=4:128] {2e-1*8/x};

%   \nextgroupplot[
%     width = 0.52\textwidth,
%     xmode = log,
%     ymode = log,
%     xmin = 1e9,
%     xmax = 1e10,
%     domain = 1e9:1e10,
%     % ymin = 1e-2,
%     % ymax = 3e3,
%     title = 128 \gls{mpi} processes,
%     xlabel = {Number of nonzero elements},
%     ylabel = {Wallclock time [\si{\second}] / Execution},
%     grid = none
%   ]
  
%   % data
%   \addplot table [y=solve, x=nonzero_elements, col sep=comma] {data/stokes/cycles_proc128.csv};

%   \addplot table [y=assemble_system, x=nonzero_elements, col sep=comma] {data/stokes/cycles_proc128.csv};
  
%   \addplot table [y=setup_system, x=nonzero_elements, col sep=comma] {data/stokes/cycles_proc128.csv};
  
%   \addplot table [y=distribute_dofs, x=nonzero_elements, col sep=comma]
%   {data/stokes/cycles_proc128.csv};
  
%   \addplot table [y=estimate_mark, x=nonzero_elements, col sep=comma] {data/stokes/cycles_proc128.csv};
  
%   \addplot table [y=refine, x=nonzero_elements, col sep=comma] {data/stokes/cycles_proc128.csv};
  
%   % auxiliary lines
%   \begin{scope}
%     % note: we always have more than 1e6 nonzeros per process !
%     \draw[dashed, very thick] ({axis cs:128000000,0}|-{rel axis cs:0,1}) -- ({axis cs:128000000,0}|-{rel axis cs:0,0});
%   \end{scope}
  
%   % optimal
%   \addplot[thick, samples=2] {1e-8*x};
%   \addplot[thick, samples=2] {5e-10*x};
%   \addplot[thick, samples=2] {1e-11*x};
  \end{groupplot}
  
  % shared legend
  \path (group c1r1.north east) -- node[above=5ex]{\ref*{legend:stokes-scaling}} (group c2r1.north west);
\end{tikzpicture}%
  \caption{\it Stokes problem.
  Left: Consecutive adaptation cycles with 128 \gls{mpi}
  processes. The dashed line again indicates $10^5$ \glspl{dof} per
  process; processes have more than this number only to the right of
  the line.
  Right: Strong scaling with a fixed problem of 15 million
  \glspl{dof}. Computations exceed $10^5$ \glspl{dof} per process
  only to the left of the dashed line.}
  \label{fig:stokes_scaling}
\end{figure}

The left panel of Fig.~\ref{fig:stokes_scaling} illustrates how run time scales with the size of the problem (here measured by the number of global \glspl{dof} $N_\text{\glspl{dof}}$) and again shows that most operations scale as one would expect given the results of the previous section.

The right panel of Fig.~\ref{fig:stokes_scaling} presents strong scaling data. As before, we get good strong scalability as long as the problem size per process is sufficiently large (to the left of the dashed line). At the same time, the figure also illustrates the limitations imposed by the linear solver we use and that have prevented us from considering larger problems: much larger problems would have taken many hours to solve even with large numbers of processes. We did not think that the associated expense in CPU cycles would have provided further insight that is not already clear from the results of the previous section and the figure -- namely, that with the exception of the linear solver and possibly assembly, all \hp-related operations scale reasonably well to large problem sizes for both simple (2d Laplace) and complex (3d Stokes) problems.

% -------------------
% --- Conclusions ---
% -------------------

\section{Conclusions}
\label{sec:conclusions}

In this manuscript, we have presented algorithms combining our
previous work on parallel and \hp-adaptive finite element methods, and
that allows us to solve problems with \hp-adaptive methods on large, parallel machines with distributed memory. In particular, we have presented algorithms for the enumeration of \glspl{dof}, a heuristic approach to weighted load balancing, and on how to transfer data of variable size between processes. 

The results we have shown in Section~\ref{sec:results} illustrate that our algorithms all scale reasonably well both to large problems and large MPI process counts, and in particular -- as one might have expected -- that (i) the linear solver is the bottleneck in solving partial differential equations that result from \hp-discretizations, and that (ii) the enumeration algorithm of Section~\ref{sec:enumeration} contributes to the overall run time in an essentially negligible way.

While the second of these conclusions makes clear that we have
succeeded in our algorithm design goals -- we have come up with an
algorithm for a task for which there was none before, and the
algorithm is fast enough to not be a bottleneck --,
our data also clearly points to future work necessary to make \hp-methods viable for more widespread use: we need more scalable iterative solvers and preconditioners, specifically ones that are better than the \gls{amg} ones we have used here. Such work would, for example, build on the \gls{gmg} ideas in \cite{mitchell2010, jomo2021}, or hybrid approaches like in \cite{fehn2020, brown2022}.
Furthermore, the literature suggests that the matrix-free approaches
of \cite{kronbichler2012, munch2022, brown2022} should be able to
overcome many of these solver limitations.

% ---------------------------------
% --- Acknowledgments & Funding ---
% ---------------------------------

\subsection*{Acknowledgments}

This paper is dedicated to the memory of William (Bill) F.~Mitchell (1955--2019), who for many years moved the \hp-finite element method along by providing high-quality implementations of the method through his \phaml{} software \cite{mitchell2002} when there were few other packages that one could play with. Equally importantly, in a monumental effort, he collected and compared the many different ways proposed in the literature in which one can drive \hp-adaptivity in practice. This work -- an extension of his work in the late 1980s comparing \h-adaptive refinement criteria \cite{mitchell1989} -- resulted in a comprehensive 2014 paper that in the end stood at 39 pages \cite{mitchell2014}, but for which the original 2011 NIST report had a full 215 pages \cite{mitchell2011}. Other significant contributions of Bill are in the area of parallelization and load balancing in form of the REFTREE algorithm \cite{mitchell2007}.

Computational methods only gain broad acceptance when the literature contains incontrovertible evidence in the form of comparison \textit{between} methods. Papers that do such comparisons are tedious to write and often not as highly regarded as ones that propose new methods, but crucial for our community to finally see which methods work and which don't. Bill excelled at writing such papers, and his contributions to \hp-finite element methods will continue to be highly regarded. His impartial and objective approach to declaring winners and losers will be missed!

An obituary for Bill Mitchell can be found at \cite{boisvert2019}.

\subsection*{Funding}

This work used compute resources provided by the Extreme Science and Engineering Discovery Environment (XSEDE) \cite{towns2014}, which is supported by National Science Foundation grant number ACI-1548562.
% see also: https://www.xsede.org/for-users/acknowledgement

MF's work was supported by the National Science Foundation under award OAC-1835673 as part of the Cyberinfrastructure for Sustained Scientific Innovation (CSSI) program.

WB's work was partially supported by the National Science Foundation under award OAC-1835673; by award DMS-1821210; and by award EAR-1925595.

% --------------------
% --- Bibliography ---
% --------------------

\bibliographystyle{ACM-Reference-Format}
\bibliography{paper.bib}

\end{document}